\documentclass[a4paper,12pt,twoside]{article}
\usepackage[latin1]{inputenc}
\usepackage[T1]{fontenc}
\usepackage{amssymb}
\usepackage{dochier}
\usepackage[final]{lpretex}
\usepackage{title}
\usepackage{a4}
\usepackage{amscd}
\usepackage[cal=boondoxo]{mathalfa}
\input xy
\xyoption{all}

\EndOfDump

\def\DHpartformat#1{(#1) }
\def\DHrefpart#1{(\DHRefpart{#1})}
\LBsetstyleof{partno}{arabic}

\def\i {\item}

\let\define\def

\def\A {{\mathbb A}}  \def\C {{\mathbb C}}

\def\GG {{\mathbb G}}   
  
  \def\P {{\mathbb P}} 
\def\Q {{\mathbb Q}} \def\R {{\mathbb R}}
\def\SS{{\mathbb S}}

\def\Z {{\mathbb Z}} 

\define \n {\mathbb N}
\define \z {\mathbb Z}
\define \q {\mathbb Q}
\define \PP {\mathbb P}

\def\sA {{\Cal A}} \def\sB {{\Cal B}} \def\sC {{\Cal C}}
 \def\sE {{\Cal E}} \def\sF {{\Cal F}}
\def\sG {{\Cal G}}  \def\sI {{\Cal I}}
\def\sJ {{\Cal J}} \def\sK {{\Cal K}} \def\sL {{\Cal L}}
\def\sM {{\Cal M}} \def\sN {{\Cal N}} \def\sO {{\Cal O}}
 
\def\sS {{\Cal S}}  
\def\sV {{\Cal V}} \def\sW {{\Cal W}} \def\sX {{\Cal X}}
\def\sY {{\Cal Y}}

\define\sdiag{\mathcal{Diag}}
\define \cN {\Cal N}
\define \cf {\Cal F}
\define \cg {\Cal G}
\define \cE {\Cal E}
\define \ce {\Cal E}
\define \cc {\Cal C}
\define \cV {\Cal V}
\define \cA {\Cal A}
\define \cK {\Cal K}
\define \cO {\Cal O}
\define \cF {\Cal F}
\define \cn {\Cal N}
\define \cI {\Cal I}
\define \sP {\Cal P}
\define \sEll {\Cal{Ell}}
\define \sJE {\Cal{JE}}
\define \sGJE {\Cal{GJE}}
\define \sHyp {\Cal H\mathcal y\mathcal p}
\def\a {\alpha}  \def\g {\gamma}  
   
\def\s {\sigma}

\define \x {\xi}
\define \y {\eta}
\define \G {\Gamma}
\define \r {\rho}
\define \w {\omega}

\def\tX {\widetilde X} 
\def \tY {\widetilde Y}

\def \tD {\widetilde D}
\def \tC {\widetilde C}
\def \tV {\widetilde V}
\def \trho {\tilde {\rho}}

\define \tH {\widetilde H}
\define \tG {\widetilde{\Gamma}}
\define \tW {\widetilde W}
\define \tF {\widetilde F}
\define \tm {\tilde m}
\define \St {\widetilde S}
\define \Xt {\widetilde X}
\define \tS {\widetilde S}
\define \tpsi {\tilde \psi}
\define \tL {\widetilde L}
\define \tE {\widetilde E}
\define \tl {\tilde l}
\define \tA {\widetilde A}
\define \tom {\tilde\omega}
\define \tT {\widetilde T}
\define \tB {\widetilde B}
\define \tf {\tilde f}
\define \tsA {\widetilde{\sA}}

\define \tM {\widetilde M}
\define \tpsi {\widetilde{\psi}}
\define \trho {\widetilde{\rho}}
\define \tR {\widetilde R}
\define \tp {\widetilde p}
\define \tq {\widetilde q}
\define \tc {\tilde c}
\define \tsF {\widetilde {\sF}}
\define \tsM {\widetilde {\sM}}
\define \tii {\tilde i}
\define \tx {\tilde x}
\define \tg {\tilde g}
\define \tw {\tilde w}
\define \tz {\tilde z}
\define \ta {\widetilde\alpha}
\define \teta{\widetilde{\eta}}
\define \tphi{\widetilde{\phi}}
\define \tsJE{\widetilde{\sJE}}
\define \ts{\widetilde{\sigma}}
\define \tsC{\widetilde{\sC}}
\define \bD {\overline{D}}
\define \bG {\overline{G}}
\define \bI {\overline{I}}
\define \bJ {\overline{J}}
\define \bK {\overline{K}}
\define\bsM{\overline{\sM}}
\define \bR {\overline{R}}
\define \bV {\overline{V}}
\define \bX {\overline{X}}
\define \bY {\overline{Y}}
\define \btau {\overline{\tau}}

\def\pd {\partial}

\def \Dx1 {\frac{\pd}{{\pd} x_1}}
\def \Dy1 {\frac{\pd}{{\pd} y_1}}
\def \Dz1 {\frac{\pd}{{\pd} z_1}}
\def \Dx2 {\frac{\pd}{{\pd} x_2}}
\def \Dy2 {\frac{\pd}{{\pd} y_2}}
\def \Dz2 {\frac{\pd}{{\pd} z_2}}

\def\q {\quad} 

\def\Mapdiagr#1{\nearrow\rlap{$\lower 5pt\vbox{{\hbox{$\mkern
-15mu\scriptstyle#1$}}}$}} 

\def\Mapdiagl#1{\llap{$\lower 5pt\vbox{{\hbox{$\scriptstyle#1\mkern
-15mu$}}}$}\searrow} 

\def\Mapswr#1{\swarrow\rlap{$\lower 5pt\vbox{{\hbox{$\mkern
-15mu\scriptstyle#1$}}}$}}              

\def\Mapnwl#1{\nwarrow\rlap{$\lower 5pt\vbox{{\hbox{$\mkern
-15mu\scriptstyle#1$}}}$}}

\def\i.e#1#2#3{\mathrel{\smash{\mathop{#2}\limits^{#1}_{#3}}}}

\def \inj {\hookrightarrow}

\define \Rhook {\hookrightarrow}

\def \half {\raise1pt\hbox{$\scriptstyle
        \frac{1}{2}\displaystyle$}} 

\def \third {\raise1pt\hbox{$\scriptstyle
    \frac{1}{3}\displaystyle$}}
\def \twothirds {\raise1pt\hbox{$\scriptstyle
    \frac{2}{3}\displaystyle$}}
\def \quarter {\raise1pt\hbox{$\scriptstyle
    \frac{1}{4}\displaystyle$}}
\def \threequarters {\raise1pt\hbox{$\scriptstyle
    \frac{3}{4}\displaystyle$}}

\def \sixth {\raise1pt\hbox{$\scriptstyle
    \frac{1}{6}\displaystyle$}}
\def \fivesixths {\raise1pt\hbox{$\scriptstyle
    \frac{5}{6}\displaystyle$}}

\def \x{{\sl X}\llap{$\mkern -2mu {\scriptstyle -}$}}

\def \End {\operatorname{End}}
\def \Hom {\operatorname{Hom}}

\def \Proj {\operatorname{Proj}}
\def \Symm {\operatorname{Sym}}
\def \Res {\operatorname{Res}}

\def \Bl {\operatorname{Bl}}

\define \Kod {\operatorname{Kod}}
\define \dimension {\operatorname{dim}}
\define \codim {\operatorname{codim}}
\define \contr {\operatorname{contr}}
\define \rk {\operatorname{rank}}
\define \Im {\operatorname {Im}}
\define \Mor {\operatorname{Mor}}
\define \Cl {\operatorname{Cl}}
\define \Hilb {\operatorname{Hilb}}
\define \degree {\operatorname{deg}}
\define \mult {\operatorname{mult}}
\define \Aut {\operatorname{Aut}}
\define \NS {\operatorname{NS}}
\define \Gal {\operatorname{Gal}}
\define \ch {\operatorname{char}}
\define \Jac {\operatorname{Jac}}
\define \Km {\operatorname{Km}}
\define \Sec {\operatorname{Sec}}
\define \Stab {\operatorname{Stab}}
\define \Br {\operatorname{Br}}
\define \Inv {\operatorname {Inv}}
\define \tr {\operatorname{tr}}
\define \Frob {\operatorname{Frob}}
\define \Symn {\operatorname{Symm}^n}
\define \Ev {\sE^\vee}
\define \ordp {\operatorname{ord}_p}
\define \Supp {\operatorname{Supp}}
\define \Ann {\operatorname{Ann}}
\define \disc {\operatorname{disc}}
\define \lie {\operatorname{lie}}
\define \embdim {\operatorname{embdim}}

\def\Def{\operatorname{Def}}

\def\Sing{\operatorname{Sing}}

\def\sEll{\sE\ell\ell}

\def\hX{\widehat{X}}

\def\bsY{\overline{\sY}}
\def\bphi{\overline{\phi}}
\def\tsY{\widetilde{\sY}}
\define\nbd{neighbourhood }

\def\hod#1#2#3#4{\ensuremath{ if#30 H^{#2}({#1},{\cal O}_{#1}) \else 
 H^{#2}(#1,\Omega^{#3}\if\relax{#4}\relax_{#1}\else _{#1/#4}\fi)\fi}}

\begin{document}
\title[Asymptotic period relations]
{Asymptotic period relations for Jacobian elliptic surfaces}
\author{N. I. Shepherd-Barron}
\address{King's College,
Strand,
London WC2R 2LS,
U.K.}
\email{Nicholas.Shepherd-Barron@kcl.ac.uk}
\maketitle
\begin{abstract}
We describe the image 
of the locus of hyperelliptic curves of genus $g$
under the period mapping in a neighbourhood
of the diagonal locus $\mathfrak{Diag}_g$.
There is just one branch for each
of the alkanes $\textrm{C}_g\textrm{H}_{2g+2}$
of elementary organic chemistry, and each branch has
a simple linear description in terms of the entries of the period matrix.

This picture is replicated for simply connected Jacobian elliptic surfaces,
which form the next simplest class of algebraic surfaces after K3 and abelian 
surfaces. In the period domain for such surfaces
of geometric genus $g$ there is a locus $\sW_{1^g}$
that is analogous to $\mathfrak{Diag}_g$,
and the image of the moduli space under the period map has just one branch
through $\sW_{1^g}$ for each alkane. Each branch is smooth and
has an explicit description
as a vector bundle of rank $g-1$ over a domain that contains $\sW_{1^g}$.
\end{abstract}
\textup{2010} \textup{MSC}:
14H42, 14J27, 32G20 (primary)
\begin{section}{Introduction}\label{intro}
The classical Schottky problem is that of describing the 
period locus $\mathfrak J_g$ of period matrices
of complex algebraic curves 
of genus $g$ as a subvariety
of Siegel space $\mathfrak H_g$.

This problem naturally extends to higher dimensions; for example,
an algebraic surface of positive geometric genus $g$
has a period matrix that arises from integrating $2$-forms around $2$-cycles,
and then the Schottky problem becomes that of describing
the image of the moduli space under the multi-valued period map.
This image we refer to as the period locus.
As described below, we are only concerned here with local aspects
of the geometry of the situation, for which the fact that the period map
is multi-valued is irrelevant.
 
We consider this problem
for simply connected Jacobian elliptic surfaces 
of geometric genus $g$. In the classification of
surfaces these are the simplest beyond 
K3 and  abelian surfaces (or those which are quotients of such surfaces);
for these latter varieties the Schottky problem scarcely arises, since
the period map is an isomorphism.
Of course, beyond those lie the surfaces of general type,
the well known complexity of whose moduli spaces
suggests that it is reasonable to focus on some
particular classes of surfaces such as the ones considered here.

The main result of this paper is that the picture for these elliptic surfaces
is analogous to that for hyperelliptic curves
and that both pictures are described by the \emph{alkanes}
of elementary organic chemistry. These are the acyclic saturated
hydrocarbons and their molecular formula is
$\textrm{C}_g\textrm{H}_{2g+2}$. The connexion with algebraic geometry is that on a hyperelliptic curve 
of genus $g$ the hyperelliptic involution has $2g+2$ fixed points
and that when such curves degenerate to trees of elliptic curves
then each elliptic curve $E$ plays the r{\^o}le of a quadrivalent carbon atom $\textrm{C}$
where the $2$-torsion points of $E$ appear as the bonds
(either $\textrm{C}$--$\textrm{C}$ or $\textrm{C}$--$\textrm{H}$) of $\textrm{C}$.

In fact,
there is a subdomain $\sW_{1^g}$ of the period domain 
$\sV_g$ for simply connected Jacobian elliptic surfaces which
corresponds to trees of $g$ \emph{special Kummer surfaces},
whose definition is recalled below, and which is isomorphic to 
$\mathfrak H_1^{g+1}$. We regard this as the analogue of the locus
$\mathfrak{Diag}_g$ of diagonal matrices in $\mathfrak H_g$, 
which corresponds to trees of elliptic curves. In both cases
there is an action of the symmetric group $\Symm_g$
on the germs $(\sV_g,\sW_{1^g})$ and $(\mathfrak H_g,\mathfrak{Diag}_g)$
that stems from the fact that the stabilizer of $\sW_{1^g}$ 
(respectively, $\mathfrak{Diag}_g$) in the relevant discrete group
is $(SL_2(\Z)\wr\Symm_g)\times SL_2(\Z)$
(respectively, $SL_2(\Z)\wr\Symm_g$),
where ${}\wr{}$ denotes the wreath product.

In slightly more detail, the main result can be summarized like this.
Let $\sJE_g$ denote the stack of simply connected Jacobian elliptic surfaces
and $PL_g$, the period locus, its image under the period map. Then
\begin{enumerate}
\item 
the branches of $PL_g$ 
through $\sW_{1^g}$ 
and the branches of the period locus $\mathfrak{Hyp}_g$
of hyperelliptic curves through $\mathfrak{Diag}_g$,
taken modulo the action of $\Symm_g$, 
both correspond to the alkanes $\textrm{C}_g\textrm{H}_{2g+2}$ and

\item each of these branches has, to first order, 
a straightforward and explicit linear description
in terms of matrices.
\end{enumerate}
In outline, the proof goes as follows: an elaboration of a plumbing
construction introduced by Fay \cite{F1}, as corrected in \cite{F2},
leads to the construction and description of one branch for each alkane,
and then a stable reduction theorem based on the Minimal Model Program (MMP)
establishes that there are no further branches.

In fact, our approach also provides a slight variant of Chakiris'
proof \cite {C1} \cite{C2} of the generic Torelli theorem for these surfaces.
\begin{notation} If $\sX$ is a separated Deligne--Mumford stack, then
we denote its geometric quotient by $[\sX]$.
if $G$ is a finite group acting on a space $X$, then
$X/G$ will denote the quotient stack and $[X/G]$
the geometric quotient.
\end{notation}
\end{section}
\begin{Acknowledgements}
I am grateful to Paolo Cascini, Bob Friedman, Mark Gross, Dave
Morrison and Richard Taylor
for some valuable conversation and correspondence
and to Hershel Farkas, Sam Grushevsky and Riccardo Salvati Manni
for their encouragement.
\end{Acknowledgements}
\begin{section}{Some further details}
We begin by recalling some constructions
that involve algebraic curves. Later we shall extend these constructions
to include algebraic surfaces.

Fay \cite{F1} 
constructed certain degenerating families
of complex algebraic curves (that is, compact Riemann surfaces)
via means of explicit plumbing constructions
that are recalled below
and then derived formulae for
the derivative of the period matrix
of each of these families.

Akira Yamada \cite{Y} then pointed out that Fay's formulae
are wrong, and gave correct formulae for Fay's constructions.

In \cite{F2} (bottom of page 123), Fay corrects his error
by pointing out that his plumbing constructions should
have been done differently, and that for these
different plumbing constructions his formulae
are correct. 

In short,
the resulting confusion can be resolved as follows.
Fay in \cite{F1} in fact made
two different plumbing constructions
of $1$-parameter degenerating families of curves
without monodromy (so the curve degenerates
but its Jacobian does not)
for which there are explicit formulae for
the derivative of the period matrix.
For one construction the correct formula is
that given in \cite{F1} and for the other
the correct formula is given in \cite{Y}. 
(There are also two different plumbing
constructions of families of curves with monodromy
but we shall not use such constructions
in this paper.) 
We shall refer to the plumbing constructions
for which Yamada's formula \cite{Y} is correct 
as Yamada plumbings,
and those for which 
Fay's formula \cite{F1}
is correct as Fay plumbings.

We shall recall the detailed construction
of Fay plumbings in Section \ref{plumbing curves}. 
In Section \ref{Fay's formulae} we will give his formula
for the derivative of the resulting period matrix;
it is convenient to point
out here that in our version of his formula there is a minus
sign that does not appear in \cite{F1}.
This is because we have chosen a different
normalization which slightly increases
the flexibility of the construction.

In fact, and this is the crux of this paper,
we do this in higher dimensions;
one advantage of Fay plumbings
is that they can be generalized to plumb not only curves
but also, at least in certain circumstances, morphisms from curves to stacks.

In the body of the paper we discuss surfaces first,
and then, in Sections \ref{curves1} to \ref{curves3}, proceed to consider curves.
We do this
because the formulae for curves are, in essence, special cases of those for surfaces.
However, in this introduction we shall take curves first.

After recovering Fay's formulae
we also recover his version \cite{F1} of
Poincar{\'e}'s ``asymptotic period relations''.
These were discovered by Poincar{\'e} \cite{P} when $g=4$
and generalized by Fay
to all values of $g$. According to Igusa's account
(see p. 167 of \cite{I}) Poincar{\'e}
exploited the geometry of the theta divisor of a Jacobian
(specifically, that it is a hypersurface of translation type),
so his argument cannot extend to the case of surfaces,
while Fay uses his plumbing construction.
We point out how these relations describe, both intrinsically and in terms
of co-ordinates, the tangent cone to the closure
$\mathfrak J_g^c$ of the Jacobian locus
along the locus $\mathfrak{Diag}_g$ of diagonal
matrices in Siegel space $\mathfrak H_g$
in terms of the Grassmannian $Grass(2,g)$ that classifies
lines in $\P^{g-1}$.

These Poincar{\'e}--Fay asymptotic period relations have also been recovered
by Farkas, Grushevsky and Salvati Manni \cite{FGSM}
in the course of proving their global weak
solution to the Schottky problem. More precisely,
they show that differentiating the identities obtained by
substituting the Schottky--Jung
proportionalities into Riemann's quartic theta identities
leads to the Poincar{\'e}--Fay relations; they deduce
their global weak solution as an immediate consequence.

Then we describe, to first order along $\mathfrak{Diag}_g$, 
the closure $\mathfrak{Hyp}_g^c$
of the hyperelliptic locus in $\mathfrak H_g$.
This description is given in terms of 
the \emph{alkanes}
of elementary organic chemistry, which were first enumerated by Cayley;
see sequence A000602 in \cite{OEIS} for corrections.
These are the hydrocarbons whose molecular formula is of the type $\textrm{C}_g\textrm{H}_{2g+2}$
(and so are exactly the saturated acyclic hydrocarbons)
and have uses ranging from fuel to furniture polish, depending on 
their molecular weight.
The $g$-alkane is the one whose carbon skeleton is a chain of length $g$.
It is distinguished from the others of the same molecular formula
by having a higher boiling point. We shall refer to $g$ as the genus of
the alkane. The number $2g+2$ is the number of fixed points of the
hyperelliptic involution on a hyperelliptic curve
of genus $g$.
It turns out that, modulo the action of the symmetric
group $\Symm_g$ on $\mathfrak H_g$
given by $s(\tau_{ij})=\tau_{s(i),s(j)}$, there is one branch of $\mathfrak{Hyp}_g^c$
along $\mathfrak{Diag}_g$ for each alkane, each branch is smooth
and there are explicit first-order equations for each branch in terms
of the entries $\tau_{ij}$ of the period matrix.

Here are more precise statements.
Recall that a square matrix $(a_{ij})$ is \emph{tridiagonal}
if $a_{ij}=0$ whenever $\vert i-j\vert\ge 2$
and \emph{quadridiagonal} if $a_{ij}=0$ whenever $\vert i-j\vert\ge 3$.
The locus of tridiagonal symmetric $g\times g$ matrices
has dimension $2g-1$.

\begin{theorem} (a special case of Theorem \ref{4.10}) To first order
(that is, modulo the square of the defining ideal)
the branch of $\mathfrak{Hyp}_g^c$ through 
$\mathfrak{Diag}_g$ in $\mathfrak H_g$
that corresponds to the $g$-alkane
is the locus of symmetric tridiagonal matrices.
\noproof
\end{theorem}

There is a similar, although slightly more complicated,
description of the branches corresponding to the other alkanes.
Poincar{\'e}'s asymptotic period relations suggest that,
on the other hand,
the fact that the locus of symmetric quadridiagonal matrices
has dimension $3g-3$ has no parallel
significance.

Now turn to surfaces.
We use Fay's plumbing to make similar
constructions and calculations for degenerating
families of 
simply connected Jacobian elliptic surfaces,
which can be regarded as the simplest surfaces 
of strictly positive Kodaira dimension and also as 
the simplest surfaces for which, thanks mainly to our understanding
of the period map for K3 surfaces, the Schottky problem
is not vacuous.

Fix an integer $h\ge 2$ and consider Jacobian elliptic surfaces 
$X$ over $\P^1$
of geometric genus $h$. These have $10h+8$ moduli,
and the coarse moduli space is rational. The primitive (co)homology
$H=H_{prim}$ of $X$ is the orthogonal complement in $H^2(X)$
of the section and a fibre; $h^2(X)=12h+10$ and $\rk H=12h+8$.
We can describe a chart of the period domain 
$$\sV_h=\left\{\xi\in Grass(h,H_{\C}):
(u,v)=0\ \textrm{and}\ (u,{\overline{u}})>0
\ \forall u,v\in\xi\right\}$$ as follows.

Pick a totally isotropic sublattice $L$ of $H$ whose rank is $h$.
\begin{definition}\label{L general} The surface $X$ is in 
\emph{$L$-general position}
if the pairing
$L\otimes\C\times H^0(X,\Omega^2_X)\to\C$
given by integration is non-degenerate.
\end{definition}
Assume that $X$ is in $L$-general position. Then a
basis $(A_1,...,A_h)$ of $L$
defines a basis $(\omega_1,...,\omega_h)$
of $H^0(X,\Omega^2_X)$ which is normalized by the requirement
$\int_{A_i}\omega_j=\delta_{ij}$. Extend $(A_1,...,A_h)$
first to a basis of $L^\perp$ and then to a basis of $H$
such that the induced basis of $H/L^\perp$ is dual to
$(A_1,...,A_h)$
and $H$ is decomposed as
$$H=L\oplus (L^{\perp}/L)\oplus(H/L^{\perp}).$$
Then the normalized period matrix of $X$ has, when we ignore the $h\times h$
identity matrix that arises from integrating around the cycles $A_i$, 
two blocks: the first
is an $h\times (10h+8)$ block that arises from integrating the forms $\omega_i$
around cycles in $L^\perp/L$
and the second is an $h\times h$ block that is skew-symmetric.
So, if $h=1$, this last block is zero and can be ignored to yield
a vector of length $18$. Different isotropic lattices $L$ 
will give different charts.

Let $PL_h$, the \emph{period locus}, denote the image of the moduli stack
in the period domain $\sV_h$ under the period map.
We shall recall the definition of the domain of the period map
later; the presence of RDPs creates a slight subtlety.

\begin{definition} A Jacobian elliptic surface is \emph{special}
if it is birational to a geometric quotient
$[C\times E/\iota]$ where $C$ is an elliptic or hyperelliptic curve,
$E$ is an elliptic curve
and $\iota$ acts on $C$ as $(-1_C)$ if $C$ is elliptic and as 
its hyperelliptic involution otherwise, and on
$E$ as $(-1_E)$.
\end{definition}

Special surfaces $X$ are characterized as the simply connected Jacobian
elliptic surfaces whose global monodromy on the cohomology
of the generic fibre
of the elliptic fibration is of order two \cite{C1}. Their geometric genus
is given by $h=p_g(X)=g(C)$. Their bad fibres are all of type $\tD_4$
(or $I_0^*$ in Kodaira's notation) and there are exactly
$2h+2$ such fibres; they correspond to the fixed points
of $\iota$.

Here are some trivial remarks, some definitions and some notation.
\begin{definition}\label{defn of sW}
\begin{enumerate}
\item $\dim\sV_h=h(10h+8)+h(h-1)/2.$

\item If $h=\sum h_i$ then
the period domain $\sV_h$ contains a copy
of $\sV_{h_1}\times\cdots\times\sV_{h_r}$.

\item Each $\sV_{h_i}$ contains a copy
$PL_{h_i,special}$ 
of the period locus
of special elliptic surfaces
of genus $h_i$.
$PL_{h_i,special}$ is isomorphic to
$\mathfrak{Hyp}_{h_i}\times\mathfrak H_1$ and so has dimension $2h_i$.

\item If $j\in [0,4]$ then $\sK_j$ denotes the period domain for Jacobian
elliptic K3 surfaces with at least $j$ fibres of type $\tD_4$.
So $\sK_j$ is a subdomain of $\sV_1$, $\dim\sK_j=18-4j$ and
$\sK_4=PL_{1,special}$.

\item Fix an alkane $\G$ of genus $h$ and let $\g_j$, for 
$j\in [1,4]$, denote the number of vertices 
(carbon atoms) in $\G$
that are joined to $j$ other carbon atoms in $\G$.
So, in particular, $h=\sum\g_j$ and $h-1=\sum j\g_j$.
There is then a closed subvariety
$\sV_\Gamma$ of the product
$\sK_1^{\g_1}\times\cdots\times\sK_4^{\g_4}$
(which is, in turn, a subdomain of $\sV_h$)
that is defined
by the requirement that the $\tD_4$-fibres on K3s
that are adjacent in $\G$
should be isomorphic.
It is easy to see that $\dim\sV_\G=9h+9$.

\item In $PL_{h_1,\ special}\times\cdots\times PL_{h_r,\ special}$,
which is isomorphic to
$\prod\mathfrak{Hyp}_{h_i}\times\mathfrak H_1^r$
and so has dimension $2h$, there is
a subvariety $\sW_{h_1,...,h_r}$ defined by the property that the
factors in each copy of $\mathfrak H_1$ are equal.
This is the period locus for unions of special elliptic surfaces
that are birational to geometric quotients $[C_i\times E_i/\iota]$
of genera $h_1,...,h_r$ where the elliptic curves $E_i$
are isomorphic, so that $\dim\sW_{h_1,...,h_r}=2h-(r-1)$.
In particular, $\dim\sW_{1^h}=h+1$.

\item For each $\G$, $\sW_{1^h}$ is isomorphic to the subvariety of $\sV_\G$
defined by the condition that each K3 surface should be special.
\end{enumerate}
\end{definition}
We regard $\sW_{1^h}$, which is isomorphic to $(\mathfrak H_1)^{h+1}$, 
as the analogue
in the period domain $\sV_h$
of the diagonal locus $\mathfrak{Diag}_g$ in
$\mathfrak H_g$.
There is an action of $(SL_2(\Z)\wr\Symm_h)\times SL_2(\Z)$ on the pair
$(\sV_h,\sW_{1^h})$,
where the action on $\sW_{1^h}$ arise from the permutation action
of $\Symm_h$
on the first term in the isomorphism $\sW_{1^h}\cong 
(\mathfrak H_1)^h\times\mathfrak H_1$.

The next result is stated in terms of a certain
vector bundle $E_\G\to\sV_\G$ of rank $h-1$. 
The fibre of $E_\Gamma$ over a point $(Y_1,...,Y_h)$ of $\sV_\Gamma$
is the vector space spanned by certain matrices
$\Pi_e$ of rank $1$, where $e=(i,j)$ runs over the edges (the carbon-carbon bonds)
of $\Gamma$. Each $\Pi_e$ is a tensor product
$\Pi_e={\underline{\omega}}_e\otimes{\underline{I}}_e$
of two vectors,
where each vector is computed from the surfaces $Y_i$
and $Y_j$. The vector ${\underline{\omega}}_e$
is comprised of projective data while 
${\underline{I}}_e$ is a vector of integrals,
so consists of transcendental data. 
These vectors are 
described explicitly in Proposition \ref{12.9}.

\begin{theorem} (= Theorem \ref{8.y}) Fix an alkane $\G$.

\part[i] There is a branch $B_\Gamma$ of $PL_h$ through $\sW_{1^h}$
that contains $\sV_\G$.

\part[i] To first order, $B_\Gamma$
is, in a neighbourhood of $\sW_{1^h}$, the vector bundle $E_\G$
over $\sV_\G$.

\part[iii] The zero section of $E_\G$ is
the image of $\sV_\G$ embedded in $B_\G$.
\noproof
\end{theorem}

This leads to the main result. It
is an analogue of Theorem \ref{4.10}.
\begin{theorem}\label{conj 2}(= Theorem \ref{main})
\part[i] The branch $B_\G$ described above is the unique branch
of the period locus $PL_h$
that contains $\sV_\G$.

\part[ii] In a neighbourhood of $\sW_{1^h}/\Symm_h$,
the period locus $PL_h$ is, to first order, 
the union of the vector bundles $E_\G$.
\noproof
\end{theorem}
\end{section}
\begin{section}{The domain of the period map for Jacobian elliptic surfaces}
  \label{the domain}
In this section we elaborate the point that, while for moduli spaces of polarized surfaces
it is natural to allow the surfaces to have RDPs, for the period map
it is better to resolve the RDPs and to allow quasi--polarizations.

We shall refer to a given Jacobian elliptic surface $f:X\to C$ with section $C_0$
as smooth if it is smooth and also relatively minimal
and as an RDP surface if it has only RDPs, $C_0$ lies in the
smooth locus of $f$ and $C_0$ is $f$-ample.
That is, if every geometric fibre of $f$ is a reduced and irreducible
curve of arithmetic genus $1$.
They are objects of two of the stacks that we shall consider:
\begin{enumerate}
\item
$\sJE^{sm}$ (resp., $\sJE^{sm}_h$)
is the stack of smooth Jacobian elliptic surfaces
(resp., simply connected such surfaces of geometric genus $h$);
\item $\sJE^{RDP}$ (resp., $\sJE^{RDP}_h$) is
the stack of RDP Jacobian elliptic surfaces
(resp., simply connected such surfaces of geometric genus $h$).
\end{enumerate}
So $\sJE^{RDP}$ is separated and there is a natural
morphism $\sJE^{sm}\to\sJE^{RDP}$ given by passing to
the relative canonical model. 
According to Artin's results \cite{Ar}, which we now recall, 
this morphism is representable,
1-to-1 on field-valued points but not separated. 

Suppose that $(X\to C\to S,C_0)$ is an object of $\sJE^{RDP}$ over $S$.
Then we have
Artin's functor $Res_{X/S}$, whose $T$-points, for an $S$-scheme $T$,
are isomorphism classes of diagrams
$$
\xymatrix{
{\widetilde{X_T}}\ar[r]^{\pi}\ar[ddr]_{F}&{X_T}\ar[r]\ar[d]&{X}\ar[d]\\
&{C_T}\ar[r]\ar[d]&{C}\ar[d]\\
&{T}\ar[r]&{S}
}$$
where $F$ is smooth, projective and relatively minimal
(this is equivalent to saying that ${\widetilde{X_T}}\to C_T$
is an object of $\sJE^{sm}$), 
and $\pi$ is projective and birational, in the sense that $\pi_*\sO=\sO$.
Then $Res_{X/S}$ is represented by a locally quasi-separated 
algebraic space $R$ over $S$, such that $R\to S$
is a bijection on all field-valued points.
The morphism $\sJE^{sm}\to\sJE^{RDP}$
described above can be localized to yield an isomorphism
$\sJE^{sm}\times_{\sJE^{RDP}}S{\buildrel{\cong}\over{\to}}R$.

Now restrict attention to simply connected
Jacobian elliptic surfaces $X$ of geometric genus $p_g=h$.
Fix a unimodular lattice $\Lambda=\Lambda_h$ of rank $12h+10$
and signature $(2h+1,10h+9)$
and elements $\s,\phi\in\Lambda$
such that $\s^2=-(h+1),\s.\phi=1$ and $\phi^2=0$. 
So $\Lambda=\Z\{\s,\phi\}\perp H$ where $H=H_h$ is unimodular,
its rank is $12h+8$
and its signature is $(2h,10h+8)$. 

Assume also that $H$ is even; 
these requirements specify $\Lambda_h$ and
$H_h$ uniquely.
(For example, $H_{prim}(X)$
is even, from the fact that $c_1(X)$ is equivalent to a multiple of 
a fibre $\phi$ and the Wu formula, which says 
that $x^2+x.c_1(X)$ is even for all classes $x\in H^2(X,\Z)$.)

 Consider the subgroup $\mathfrak G$
of the orthogonal group $O_\Lambda(\Z)$ given by
$$\mathfrak G=\{\g\in O_\Lambda(\Z)\vert\g(\s)=\s\ \textrm{and}\ \g(\phi)=\phi\};$$
then $\mathfrak G$ is naturally identified with $O_H(\Z)$.
There is a $\mathfrak G$-torsor $\sJE_\Lambda\to \sJE_h^{sm}$ defined by the fact
that a $T$-point of $\sJE_\Lambda$ consists of a $T$-point of $\sJE_h^{sm}$
and an isometry $\Psi:\Lambda_T\to R^2F_*\Z$ such that
$\Psi$ maps $\s$ to the class of the section $(C_0)_T$
of ${\widetilde{X_T}}\to C_T$
and $\phi$ to the class of a fibre.

The discussion in the third paragraph on p. 228 of \cite{C1}
can be translated into
the language of stacks to say that
$\sJE_\Lambda$ is the domain of the period map.
That is, the period map is a $\mathfrak G$-equivariant
holomorphic morphism
$$\widetilde{per}:\sJE_\Lambda\to\sV_h,$$
where $\sV_h$ is the period domain.
Equivalently, the period map is the quotient morphism
$$per:\sJE^{sm}_h=\sJE_\Lambda/\mathfrak G\to\sV_h/\mathfrak G$$
of quotient stacks. This fits into a commutative diagram
$$\xymatrix{
{\sJE_h^{sm}}\ar[r]^{per}\ar[d]&{\sV_h/\mathfrak G}\ar[d]\\
{\sJE^{RDP}_h}\ar[r]&{[\sV_h/\mathfrak G].}
}$$
If $S$ is the henselization of $\sJE^{RDP}_h$ at a closed point,
so that $S$ is the base of a miniversal deformation 
of a simply connected RDP Jacobian elliptic surface $X_s$, then \cite{Ar}
the henselization $R^{hens}$ of $R=\sJE^{sm}_h\times_{\sJE^{RDP}_h}S$
at its unique closed point is the base of a miniversal
deformation of the minimal resolution $\tX_s$ of $X_s$.
(In fact Artin shows that $R^{hens}$ is the base of a versal deformation,
but in characteristic zero his argument proves the miniversality
of $R^{hens}$.)
So, for local purposes such as those of this paper,
we can take the domain of the period map to be
a miniversal deformation space of a smooth Jacobian elliptic surface.

In fact, slightly more is true. The finite Weyl group $W$
associated to the configuration of singularities on $X_s$
acts on $R^{hens}$ in such a way that $[R^{hens}/W]=S$
and there is a commutative diagram
$$\xymatrix{
{R^{hens}}\ar[r]\ar[dr]\ar@/_2pc/[ddr]&{R}\ar[dr]^{per}\ar@/_3pc/[dd]\\
&{R^{hens}/W}\ar[r]\ar[d]&{\sV_h/\mathfrak G}\ar[d]\\
&{S}\ar[r]&{[\sV_h/\mathfrak G].}
}$$
This is because 
there is a non-separated union $\tR_P=\cup_{\sC}R^{hens}_\sC$
of copies of $R^{hens}$ on which 
$W$ acts freely, where there is
one copy $R^{hens}_\sC$ of $R^{hens}$
for each chamber $\sC$ that is defined in the usual way
and two charts $R^{hens}_\sC, R^{hens}_\sE$
in $\tR_P$ are glued in a way that depends upon
the relative position of the chambers $\sC,\sE$
with respect to $W$.
(This glueing always
induces an isomorphism over the complement 
of the discriminant locus in $S$.)
Moreover $W$ permutes the charts
$R^{hens}_\sC$ simply transitively, 
and $\tR_P/W=\tsJE_h\times_{\sJE_h}S=R$. 
The isomorphisms $R^{hens}_\sC\to R^{hens}$ glue to
a morphism $\a:\tR_P\to R^{hens}$ that exhibits $R^{hens}$
as the maximal separated quotient of $\tR_P$
in the category of schemes. So
$W$ also acts on $R^{hens}$, $\a$ is $W$-equivariant and $S=[R^{hens}/W]$.
All this is documented in \cite{SB}, to which we add one comment:
since, in the commutative diagram
$$
\xymatrix{
{\tR_P}\ar[r]^{\a}\ar[d]&{R^{hens}}\ar[d]\ar[dr]&\\
{R}\ar[r]&{R^{hens}/W}\ar[r]&{S}
}
$$
the vertical maps are quotients by $W$,
so that the square is $2$-Cartesian,
it follows
that the scheme $S$ is the maximal separated quotient of 
the algebraic space $R$
in the category of algebraic spaces
and that the stack $R^{hens}/W$ is the maximal separated
quotient of $R$ in the $2$-category of algebraic stacks.

Finally, consider the forgetful morphism 
$\a:S\to\prod_x Def_{X_s,x}$
from $S$ to the product
of the miniversal deformation spaces of the 
germs $(X_s,x)$ of $X_s$ at its singularities $x$
and the morphisms $Def_{\tX_s,x}\to Def_{X_s,x}$.
These latter are geometric quotients by the relevant finite Weyl group
and $R^{hens}$ is isomorphic to the fibre product
$$S\times_{\prod_x Def_{X_s,x}}\prod_x Def_{\tX_s,x},$$
so that
if $\a$ is smooth then so is $R^h$.
\end{section}
\begin{section}{Fay's plumbing for curves and stacky curves}\label{plumbing curves}
In this section we give a detailed description of Fay
plumbings. All of this is taken from \cite{F1} and \cite{F2} but we repeat it here
to avoid confusion.

Fix a real number $\delta$ with $0<\delta\ll 1$.
Let $\Delta=\Delta_t$ be the open disc in the complex plane $\C$
with co-ordinate $t$ defined by $\vert t\vert<\delta^2$
and let $F$ be the open submanifold of $\C^2$
with co-ordinates $q,v$ defined by 
$\vert q\vert<\delta^{1/2}$, $\vert v\vert<\delta$
and $\vert v\pm q\vert<\delta$
(so that, in particular, $\vert v^2-q^2\vert <\delta^2$).
Set $t=v^2-q^2$. Note that the
morphism $t:F\to\Delta$ is smooth outside the origin in $\Delta$
and the fibre $t^{-1}(0)$ consists of two discs
$U_a'$ and $U_b'$, where $U_a'$ is given by $q-v=0$
and $U_b'$ by $q+v=0$. These discs cross normally at the
point $0$ given by $q=v=0$.

We shall use $F$, with these co-ordinates, 
as the basic plumbing fixture.

Now suppose that $C_a,C_b$ are curves of genus $g_a,g_b$,
that $a\in C_a, b\in C_b$ and that if $i=a,b$ then 
we have chosen a local
co-ordinate $z_i$ on some simply connected neighbourhood $U_i$
of $i$ in $C_i$ such that $U_i$ is defined by the inequality 
$\vert z_i\vert<\delta^{1/2}$.
In particular, $U_i$ is an open disc that is embedded in $C_i$.

The suffix $i$ will stand for either $a$ or $b$ in what follows.

Let $U_i^c$ denote the closure of $U_i$ in $C_i$.
We shall assume that each $U_i^c$ is a (closed) disc;
this is slightly stronger than assuming $U_i$ is an open disc.
We shall assume also \emph{either} that $C_a$ and $C_b$ are distinct
\emph{or} that $C_a=C_b$ and 
the closed discs $U_a^c$ and $U_b^c$ are disjoint.

Of course, all of these assumptions can be fulfilled
after decreasing $\delta$ if necessary.

In $F$ we have open subsets $V_a$ and $V_b$
defined, respectively, by the inequalities
\begin{eqnarray}\label{plumbing formulae}
\vert v-q\vert&<&\vert q\vert,\ \vert t\vert<\delta \vert q\vert^2\ \textrm{and}\\
\vert v+q\vert&<&\vert q\vert,\ \vert t\vert<\delta \vert q\vert^2.
\end{eqnarray}
Note that $V_a$ and $V_b$ are disjoint
and that $V_i$ contains $U_i'-\{0\}$.

Define $\rho_i:V_i\to U_i\times\Delta$ by
$z_i=q,\ t=q^2-v^2$.
It is easy to check that $\rho_i$ is unramified
(since the ramification locus is defined by $v=0$)
and is injective. Therefore $\rho_i$ is an isomorphism from
$V_i$ to its image $Y_i$, which is open in $U_i\times\Delta$
and so open in $C_i\times\Delta$.

Let $Y_i$ denote the image of $V_i$.
Then $Y_i$ is contained in the
region $\vert t\vert<\delta\vert z_i\vert^2$
and contains $(U_i-\{i\})\times\{0\}$.
Note that $\rho_i$ maps $U_i'-\{0\}$ isomorphically onto $U_i-\{i\}$.

In $C_i\times\Delta$ consider the closed subset
$\{(P,t)\ :\ P\in U_i^c\ \textrm{and}\ (P,t)\not\in Y_i\}$.
Define $W_i$ to be the open subset of $C_i\times\Delta$
obtained by deleting this closed set.

Note that $W_i\cap (C_i\times\{0\})=C_i-\{i\}$.

\begin{proposition}
We can glue the three charts $W_a$, $W_b$ and $F$ 
via the isomorphisms $\rho_i:V_i\stackrel{\cong}{\to} Y_i$
to get a separated $2$-dimensional 
complex manifold $\sC$ with a
holomorphic map $\pi:\sC\to\Delta_t$
such that $\pi^{-1}(0)=C_a\cup_{a\sim b} C_b$,
where $C_a$ and $C_b$ cross normally at a
single point. After shrinking $\Delta_t$ if necessary,
the map is proper.
\begin{proof}
This is routine.
\end{proof}
\end{proposition}

We call this a \emph{Fay plumbing}
because it is, implicitly, constructed and considered by Fay
in the final paragraph of p. 37 of \cite{F1}.

If $C_a$ and $C_b$ are distinct then
for all $t\ne 0$ the fibre $\sC_t$
is a curve of genus $g=g_a+g_b$
and there is no monodromy on its cohomology.
This is what Fay calls ``pinching a cycle homologous to zero''
\cite{F1}, p. 37 \emph{et seq.}
We shall call it a 
\emph{homologically trivial Fay plumbing of $C_a$ to $C_b$ that identifies $a$ with $b$}
or just a \emph{Fay plumbing of $C_a$ to $C_b$} without explicit mention of the points
or co-ordinates that are chosen.

If $C_a=C_b$, of genus $g$, and also $U_a^c\cap U_b^c=\emptyset$
(so that the plumbing is possible),
then for all $t\ne 0$ the fibre $\sC_t$ is of genus $g+1$,
$\sC_0$ is the nodal curve $C/a\sim b$ and there is non-trivial
monodromy on the cohomology of $\sC_t$.
Fay calls this ``pinching a non-zero homology cycle''
\cite{F1}, p. 50.

Note that on $F\cap W_i$ we have $z_i=q$.

In a homologically trivial Fay plumbing
the derivative at $t=0$ of the period matrix
of $\sC_t$ is given as Corollary 3.2 on p.\ 41 of \cite{F1}.
We shall prove this, or, rather, a version
of it in a higher-dimensional context,
later on, in Proposition \ref{8.2}.

When a non-zero homology cycle is pinched
in a Fay plumbing the period matrix is
described by by Corollary 3.8 on p.\ 53 of \cite{F1}.

However, the families that are written down on
pp.\ 37 and 50 of \emph{loc.\ cit.},
which we call the \emph{Yamada plumbing}
because they are considered explicitly by Yamada,
are given by the glueing $t=z_az_b$ and 
(provided that $g_a+g_b\ge 1$)
are different families
because, for example, their period matrices are different.
Their expansions at $t=0$ are given as Corollary 2 on p. 129
and Corollary 6 on pp. 137-138 of \cite{Y}.

In this paper we shan't have any need to pinch
cycles that are not homologous to zero.
For one application of such pinchings, however, see \cite{CSB}.
Note that there Yamada plumbings are used,
although an argument could also be based on the use
of Fay plumbings.

It will also be useful to be able to plumb certain stacky curves.

\begin{definition} A \emph{$\Z/2$-curve} is a reduced connected proper
Deligne--Mumford stack $\tC$ of dimension $1$ such that
at each generic point the stabilizer (the isotropy group)
is trivial and there are finitely many
points where the stabilizer is $\Z/2$.
\end{definition}
For example,
if $i\in C_i$
then there is a unique smooth $\Z/2$-curve
$\tC_i$ such that the
geometric quotient $[\tC_i]$ of $\tC_i$ is given by
$[\tC_i]=C_i$, the quotient map
$\tC_i\to C_i$ is an isomorphism over $C_i-\{i\}$ and
there is a unique geometric point on $\tC_i$
lying over $i$.
There are local co-ordinates $\tz_i$ and $z_i$
on $\tC_i$ and on $C_i$, respectively, such that
the non-trivial element
$\iota$ of $\Z/2$ acts via $\iota^*\tz_i=-\tz_i$ and
$z_i=\tz_i^2$.

Take the plumbing fixture $F$, with co-ordinates $q,v$
as before, and let $\iota$ act on $F$ by $\iota^*q=-q$,
$\iota^*v=-v$. So the quotient stack $\tF=F/\langle\iota\rangle$
is a smooth separated 2-dimensional Deligne--Mumford stack.

Suppose $\tC_i$, for $i=a,b$, are smooth $\Z/2$-curves
and that each of $a,b$ the stabilizer is 
$\langle\iota\rangle\cong\Z/2$.
Fix a local co-ordinate $\tz_i$ on $\tC_i$
at $i$ such that $\iota^*\tz_i=-\tz_i$.
\begin{proposition}\label{stacky plumb}
We can glue together the charts $\tC_a$, $\tC_b$ and $\tF$
via the formulae $\tz_i=q$ and $t=q^2-v^2$ 
to get a smooth $2$-dimensional Deligne--Mumford stack $\tsC$
with a proper separated morphism $\tsC\to\Delta_t$.
\begin{proof}
As before.
\end{proof}
\end{proposition}

Taking geometric quotients 
gives the following result.
Define $G=[F/\langle\iota\rangle]$.
\begin{proposition}\label{plumbing involutions}
The curves $C_i$ can be plumbed via the plumbing fixture $G$
to $\sB\to\Delta_t$. On the surface
$\sB$ there is an $A_1$-singularity at $0$
and $\sB_0=C_a\cup C_b/a\sim b$.
\begin{proof}
This is clear.
\end{proof}
\end{proposition}
The minimal resolution of $\sB$
gives a semi-stable family of curves over $\Delta_t$.
\end{section}
\begin{section}{Plumbing families of curves}
The plumbing construction of the previous section can be extended
so as to plumb families of curves.

So suppose that $C_a\to S_a$, $C_b\to S_b$ are analytic families 
of semi-stable curves of genera $g_a$ and $g_b$,
respectively, with sections $\s_i:S_i\to C_i$
each of which lies in the relevant
smooth locus. Assume also that $U_i$ is a tubular neighbourhood
of the section $i=\s_i(S_i)$ and that $z_i$ is a fibre co-ordinate
on $U_i$ such that $U_i$ is isomorphic to 
$S_i\times\Delta_{z_i}$, where $\Delta_{z_i}$ is the
disc with co-ordinate $z_i$ such that $\vert z_i\vert<\delta$.
Such a neighbourhood and such a co-ordinate will exist 
if each $S_i$ is sufficiently small,
for example a small polydisc.

Take the same $2$-dimensional plumbing fixture $F$ 
and disjoint open subsets $V_i$ as before
and consider the morphisms 
$\rho_i:V_i\times S_i\to U_i\times\Delta_t$
defined by
$$\rho_i(q,v,s_i)=(\s_i(s_i),z_i,t)$$
where $z_i=q$ and $t=q^2-v^2$.

As before, $\rho_i$ is an isomorphism to an open analytic subvariety
$Y_i$ of $U_i\times\Delta_t$.

In $C_i\times\Delta_t$ consider the closed subset 
that is the intersection of $U_i^c\times\Delta_t$
and the complement of $Y_i$,
and define $W_i$ to be the open subvariety of
$C_i\times\Delta_t$ obtained by deleting this closed subset.

Now glue the chart $W_a\times S_b$ to $F\times S_a\times S_b$
via the isomorphism $V_a\times S_a\times S_b\to Y_a\times S_b$,
and glue $W_b\times S_a$ to $F\times S_a\times S_b$
via the isomorphism $V_b\times S_a\times S_b\to Y_b\times S_a$.

The result of this glueing is an analytic space $\sC$
with a proper morphism $\sC\to S_a\times S_b\times\Delta_t$
that is a family of semi-stable curves of genus $g_a+g_b$.
If $D_i\subset S_i$ is the discriminant locus of $C_i\to S_i$
then the discriminant locus of $\sC\to S_a\times S_b\times\Delta_t$
is $D_a\times S_b\times\Delta_t\cup D_b\times S_a\times\Delta_t
\cup S_a\times S_b\times\{0\}$.

We can concatenate plumbings as follows: 
suppose that $C_i\to S_i$ are families of stable curves
for $i=a,b,c$ and that $\s_a, \s_{b_1},\s_{b_2},\s_c$
are sections of $C_a,C_b,C_b,C_c$ respectively
and that $\s_{b_1}$ and $\s_{b_2}$ are disjoint.

Then there are two choices: we can first plumb $C_a$ to
$C_b$, obtaining $\sC'\to S_a\times S_b\times\Delta$, 
in a way that identifies the sections $\s_a$ and $\s_{b_1}$,
and then plumb $\sC'$ to $C_c$ to get a family over
$S_a\times S_b\times\Delta\times S_c\times\Delta$,
or we can first plumb $C_b$ to $C_c$, obtaining $\sC''$, by identifying
$\s_{b_2}$ with $\s_c$, and then plumb $\sC''$ to $C_a$
to get another family over $S_a\times S_b\times\Delta\times S_c\times\Delta$.
It will be important for us to notice that
these two families are the same; that is,
the final result is independent of the order of the plumbings.

Similarly we can construct stacks $\tC_i\to S_i$ by introducing
isotropy groups $\Z/2$ along each section $\s(S_i)$
and then plumb together the stacks $\tC_i$.
\end{section}
\begin{section}{Plumbing morphisms}\label{plumbing morphisms}
Keep the notation of Section \ref{plumbing curves}
and suppose that $\sC=W_a\cup W_b\cup F\to\Delta$ is a
Fay plumbing
of curves $C_a$ to $C_b$ that identifies $a$ with $b$
and that $z_i$ is the local co-ordinate on $C_i$ at $i$
that is used in the plumbing. 

Suppose that $\sM$ is some holomorphic stack
and that $\phi_i:C_i\to\sM$ are morphisms
such that for the choice of co-ordinates $z_a$ and $z_b$
the morphisms $\phi_a$ and $\phi_b$ are isomorphic on $U_a$ and on $U_b$.
That is, there is a $2$-commutative diagram
$$
\xymatrix{
{U_a}\ar[rr]^{\a}_{\cong}\ar[dr]_{\phi_a\vert_{U_a}}&&{U_b}\ar[dl]^{\phi_b\vert_{U_b}}\\
&{\sM}
}
$$
where $\a$ is an isomorphism such that $\a^*z_b=z_a$.

\begin{lemma} There is a morphism $\Phi:\sC\to\sM$
such that the restriction of $\Phi$ 
to the closed fibre $\sC_0=C_a\cup C_b$
coincides with $\phi_i$ on each $C_i$
and the restriction $\Phi_i$ of $\Phi$ to each of the two charts $W_i$
of $\sC$
factors through the projection $pr_i:W_i\to C_i$
as $\Phi_i=\phi_i\circ pr_i$.
\begin{proof} Suppose that $w_1,...,w_n$ are local co-ordinates
on a smooth chart $X\to\sM$ and that
a local lift $\tphi_i:U_i\to X$ of
$\phi_i\vert_{U_i}$ is given, in local co-ordinates,
by a formula 
$$\tphi_i(z_i)=(\phi_{i,1}(z_i),...,\phi_{i,n}(z_i)).$$
Then the hypothesis is that $\tphi_{a}(z_a)$ is equivalent
to $\tphi_{b}(z_b)$,
so that we can define $\Psi:F\to\sM$ by
$\Psi(q,v)=(\phi_{a,1}(q),...,\phi_{a,n}(q))
\sim(\phi_{b,1}(q),...,\phi_{b,n}(q)).$

Then define $\Phi_i:W_i\to\sM$ as
$\Phi_i=\phi_i\circ pr_i$ where
$pr_i:W_i\to C_i$ is the projection.
Since the morphisms $\Psi$ and $\Phi_i$
agree on the overlaps the lemma is proved.
\end{proof}
\end{lemma}
\begin{remark} If instead we consider the Yamada plumbing defined
by the glueing $t=z_az_b$ then it is not clear how to construct
such a morphism $\Phi$. This is why we focus here on Fay plumbings.
\end{remark} 
\begin{definition}\label{ell stacks}
\begin{enumerate}
\item $\sEll^0$ denotes the stack of (smooth) elliptic curves.
\item $\sEll$ is the stack of reduced and irreducible curves of genus $1$
with planar singularities, provided with a marked smooth point.
\item $\sEll^{Kod}$ is the stack of Kodaira fibres.
\end{enumerate}
\end{definition}
For example, suppose that $\sM$ is one of the stacks $\sEll$ or $\sEll^{Kod}$.
Then, for either of these choices of $\sM$, the assumptions hold
if the points $\phi_i(i)$
are both isomorphic to the same smooth elliptic curve
and each $\phi_i$ has the same finite order of ramification at
the point $i\in C_i$.

In general, the hypotheses imply that the curves $C_a,C_b$
have the same image in the geometric quotient of $\sM$
if that quotient exists.

By definition, the datum of a \emph{Jacobian elliptic fibration}
$f:X\to Y$ over an irreducible base is equivalent to the datum of a morphism
$\phi:Y\to \sEll$ such that $Y$ meets $\sEll^0$.
Then the next corollary, which is an immediate consequence of the lemma,
shows that Jacobian elliptic surfaces can be plumbed.

\begin{corollary}\label{5.2} Suppose that, for each $i=a,b$, 
$f_i:X_i\to C_i$ is a smooth (resp., RDP) Jacobian
elliptic surface defined by morphisms $\phi_i:C_i\to\sEll^{Kod}$
(resp., $\phi_i:C_i\to\sEll$),
that the curves $\phi_a(a)$ and $\phi_b(b)$ are isomorphic 
to the same smooth elliptic curve $E$ and that
the orders of ramification of
$\phi_a$ at $a$ 
and of $\phi_b$ at $b$ are finite and equal.

Then there is a Jacobian elliptic fibration
$f:\sX\to\sC$ and a proper morphism $g:\sC\to\Delta$ with
composite $h=g\circ f:\sX\to\Delta$ where
\begin{enumerate} 
\item $\sX$ is a smooth (resp., c-DV) threefold,

\item $\sC$ is a smooth surface,

\item if each $X_i$ is simply connected and $C_a\ne C_b$ 
then $g$ is a homologically trivial Fay plumbing
of two copies of $\P^1$ and so 
is a conic bundle whose singular fibre $\sC_0$ has two components,
 
\item if $C_a\ne C_b$ then the fibre $\sX_0=X_a\cup X_b=h^{-1}(0)$
is a reduced divisor in $\sX$
with normal crossings along the common curve $E$,

\item if $\sW_i$ is the inverse image in $\sX$
of the open subset $W_i$ of $\sC$,
then there is a holomorphic projection $\sW_i\to X_i-E$
such that the induced morphism $\sW_i\to (X_i-E)\times\Delta$
is an isomorphism to an open submanifold that contains
$(X_i-E)\times\{0\}$ and

\item $\sX_t\to\sC_t$ is a smooth (resp., RDP) Jacobian elliptic surface
for $t\ne 0$ and

\item the configuration of $(-2)$-curves on $\sX_t$ is constant
and specializes to the union of the $(-2)$-configurations
on $X_a$ and on $X_b$ if $X_a$ and $X_b$ are smooth,
while if they are RDP surfaces then the analogous statement
for configurations of RDPs is true.
\end{enumerate}
\begin{proof} The hypotheses concerning ramification imply
that $\phi_a$ and $\phi_b$ are locally isomorphic.
Then everything except the last statement, which is easy, is automatic.
\end{proof}
\end{corollary}
In particular, suppose that each $X_i$ is simply connected.
Then $C_i$ is isomorphic to $\P^1$
and $\sC_t$ is also isomorphic to $\P^1$.
It follows that if $\sX\to\Delta$ is the composite morphism
then for $t\ne 0$ the fibre
$\sX_t$ is a simply connected Jacobian elliptic surface
and that, if $p_g(X_i)=g_i$,
$\sX_0=X_a\cup X_b$, $\sX_0$ is reduced with normal crossings
and $X_a\cap X_b=E$. Moreover,
$p_g(\sX_t)=g_a+g_b+1=g$, say.

\begin{definition}\label{RDP,sst}
$\sJE^{RDP,sst}$ denotes the stack of 
\emph{semi-stable RDP Jacobian elliptic surfaces}.
Its objects over a scheme $B$ are pairs
$(\sX{\buildrel{f}\over{\to}}S\to B, \tS)$ such that
\begin{enumerate}
\item $f:\sX\to S$ and $S\to B$ are projective, flat, Cohen--Macaulay
and of relative dimensions $2$ and $1$, respectively,
\item $\tS$ is a section of $f$ lying in the relatively smooth locus of $f$
(so is an effective Cartier divisor on $\sX$),
\item \label{req 3} $\tS$ is $f$-ample,
\item $\sX\to B$ and $S\to B$ are semi-stable, in the sense that
all geometric fibres are reduced with normal crossings, and 
the family of curves $S\to B$
is of compact type,
\item for every geometric point $\delta$ of $B$ and every
irreducible component $C$ of $S_\delta$, $f^{-1}(C)\to C$
is an RDP Jacobian elliptic surface whose origin is
$\tS\cap f^{-1}(C)$ and
\item $f:\sX\to S$ is smooth over the nodes of $S_\delta$.
\end{enumerate}
\end{definition}
\begin{definition}\label{sm,sst}
$\sJE^{sm, sst}$ is the stack obtained from
$\sJE^{RDP,sst}$ by dropping
the requirement \ref{req 3} above 
and imposing instead the requirement
that each $f^{-1}(C)\to C$ should be a smooth and relatively
minimal Jacobian elliptic surface. 
\end{definition}
Note that $\sJE^{sm,sst}$ (resp., $\sJE^{RDP,sst}$)
contains $\sJE^{sm}$ (resp., $\sJE^{RDP}$) as an open substack.

In terms of this notation,
the objects $\sX\to\sC\to\Delta$ of Corollary \ref{5.2}
are objects over $\Delta$ of either $\sJE^{sm,sst}$ or $\sJE^{RDP,sst}$,
as appropriate.

There is an obvious extension of this 
to plumbings of several curves, as follows.

Suppose that $C_1,...,C_r$ are curves, that $a_1\in C_1,\ b_1,a_2\in C_2,\ldots,
b_{r-2},a_{r-1}\in C_{r-1}$, $b_{r-1}\in C_r$
and that $z_x$ is a local co-ordinate at $x$ on the relevant curve.
Then there is a Fay plumbing $\sC\to S_{r-1}$ where 
$S_{r-1}$ is an $(r-1)$-dimensional polydisc and the curves
$C_1,...,C_r$ are attached in a chain that identifies $a_i$ with $b_i$.
As a manifold, $\sC$ is the union of charts $W_i$, each
of which is open in $C_i\times S_{r-1}$
and contains $C_i-\{a_i,b_{i-1}\}$, together with plumbing fixtures.

Suppose also that there are given morphisms $\phi_i:C_i\to\sM$
such that $\phi_i$ and $\phi_{i+1}$ are locally isomorphic
in terms of the given co-ordinates $z_{a_i}$ and $z_{b_i}$.
Then these morphisms can be plumbed to give a morphism
$\Phi:\sC\to\sM$ such that the restriction $\Phi_i$
of $\Phi$ to each chart $W_i$
factors through the projection $pr_i:W_i\to C_i$
as $\Phi_i=\phi_i\circ pr_i$.

Moreover, these plumbings can be constructed in any order.
\end{section}
\begin{section}{Plumbings modulo $t^2$}\label{first order}
Suppose that $C_i$, for $i=a,b$, are curves,
that $i\in C_i$ and that $z_i$ is a local co-ordinate on
$C_i$ at $i$, as before, and that $\sC\to\Delta$ is the corresponding
Fay plumbing of $C_a$ to $C_b$. Suppose also that $\phi_i:C_i\to\sM$
are morphisms that are merely 
\emph{isomorphic to first order}. That is,
$\phi_a':\Sp\C[[z_a]]/(z_a^2)\to\sEll$ is isomorphic to
$\phi_b':\Sp\C[[z_b]]/(z_b^2)\to\sEll$
in the sense that there is a $2$-commutative diagram
$$\xymatrix{
{\Sp\C[[z_a]]/(z_a^2)}\ar[rr]^{\a}_{\cong}\ar[dr]_{\phi_a'}&&
{\Sp\C[[z_b]]/(z_b^2)}\ar[dl]^{\phi_b'}\\          
&{\sM}
}
$$
where $\a$ is an isomorphism such that $\a^*z_b=z_a$.
Put $\Delta'=\Sp\C[[t]]/(t^2)$ and $\sC'=\sC\times_\Delta\Delta'$.

\begin{proposition}
$\phi_a'$ and $\phi_b'$ can be plumbed to give a morphism
$\Phi':\sC'\to\sM$.
\begin{proof} Exactly as before.
\end{proof}
\end{proposition}

\begin{corollary}\label{whatnot}\label{7.2}
If $f_i:X_i\to C_i$ are Jacobian elliptic
surfaces such that the classifying morphisms 
$\phi_i:C_i\to\sEll$ are isomorphic
to first order at $a$ and $b$, then $f_a$ and $f_b$ can be plumbed to
$\sX'{\buildrel{f}\over{\to}}\sC'\to\Delta'$. 
\noproof
\end{corollary}

If each $f_i:X_i\to C_i$ is a smooth (resp, RDP) 
Jacobian elliptic surface
then $\sX'{\buildrel{f}\over{\to}}\sC'\to\Delta'$ 
is an object over $\Delta'$
of $\sJE^{sm,sst}$ (resp., $\sJE^{RDP,sst}$).

\begin{lemma} Suppose that each $f_i:X_i\to C_i$
is an RDP Jacobian elliptic surface with minimal resolution
$\tX_i\to X_i$. Then $\sX'{\buildrel{f}\over{\to}}\sC'\to\Delta'$
can be lifted to an object of $\sJE^{sm,sst}$ over $\Delta'$.
\begin{proof} By construction, $\sX'\to\Delta'$ induces a trivial
deformation of each RDP.
\end{proof}
\end{lemma}
\begin{lemma}\label{sm via sextic} 
The forgetful morphism $\sJE^{RDP,sst}\to\sM^c=\coprod\sM_g^c$
to the stack of stable curves of compact type
is smooth and has irreducible fibres.
\begin{proof}
Given an object $C$ of $\sM^c$ and an object $p:Y\to C$ of $\sJE^{RDP,sst}$
over $C$, define $\sL=p_*\sO_Y(C_0)$. This is an $\iota$-linearized
$p$-ample
invertible sheaf on $C$ and $\GG_m$ acts $\iota$-equivariantly
on the vector bundle $\sL\oplus\sL^{\otimes 2}\oplus\sL^{\otimes 3}$ with weight
$n$ on $\sL^{\otimes n}$. Define
$\P=[(\sL\oplus\sL^{\otimes 2}\oplus\sL^{\otimes 3})-\{0\}]/\GG_m$;
this is a $\P(1,2,3)$-bundle over $C$ and
$Y$ is a sextic Cartier divisor in $\P$.

Therefore an object of $\sJE^{RDP,sst}$ lying over an object $C$ of $\sM^c$
is defined by the data of a line bundle $\sL$ on $C$ and a sextic divisor
in the $\P(1,2,3)$ bundle associated to $\sL$; it follows that
$\sJE^{RDP,sst}\to \sM^c$
is smooth, with irreducible fibres.
\end{proof}
\end{lemma}
\begin{corollary}\label{construct} 
The object
$\sX'{\buildrel{f'}\over{\to}}\sC'\to\Delta'$ of $\sJE^{RDP,sst}$
over $\Delta'$ can be extended to
an object $\sX{\stackrel{f}{\to}}\sC\to\Delta$ of $\sJE^{RDP,sst}$
over $\Delta$ where the $3$-fold $\sX$ is smooth.
\begin{proof} The existence of $\sX$ is clear and
the only thing to check is that $\sX$ is smooth.
For this, notice that in suitable local co-ordinates $\sX'\to\Delta'$
is defined by $t=xy$; then the same holds for $\sX$
and the smoothness follows.
\end{proof}
\end{corollary}
In Lemma \ref{sm is smooth}
we shall prove a version of this result
for another relevant partial compactification of $\sJE^{sm}$,
which will be denoted by $\sJE^{sm,n.m.}$. The distinction between these
compactifications is that a typical boundary point
of $\sJE^{sm,sst}$ corresponds to a singular surface
that consists of two smooth surfaces that meet transversely
in a smooth elliptic curve, while in $\sJE^{sm,n.m.}$
the two surfaces are singular and meet in a model of a
$\tD_4$-fibre.
\end{section}
\begin{section}{Homologically trivial plumbings of surfaces}\label{11}
Here we construct families of Jacobian elliptic surfaces 
that degenerate in a homologically trivial fashion,
so that the period matrix specializes to a matrix
lying in the interior of the period domain
and there is an explicit formula for the derivative
of the period matrix that involves no derivatives
but is instead a matrix of rank one, just as for curves.

This leads to a result analogous to the 
asymptotic relations that will be derived later, in Section \ref{curves3}, 
for the locus of hyperelliptic curves.

We begin with a modification of Kodaira's degenerate fibres
and his notation for them.
\begin{enumerate}
\item $I_{n-4}^*=\tD_n$ for $n\ge 4$:
contract the four $(-2)$ curves at the extremities
and call the result ${\overline D}_n$. This has four $A_1$ singularities.

\item $II=Cu$: blow up the cusp three times, then contract
the resulting $(-2),(-3)$ and $(-6)$ curves, and call the result
$\overline{II}$. This has three singularities, one of each type
$\frac{1}{2}(1,1)=A_1,\frac{1}{3}(1,1)$ and $\frac{1}{6}(1,1)$.

\item $III=Ta$: blow up the tacnode twice,
then contract the resulting $(-2)$ curve and two $(-4)$ curves
and call the result $\overline{III}$. This has one
$A_1$ singularity and two singularities
of type $\frac{1}{4}(1,1)$.

\item $IV=Tr$: blow up the triple point,
then contract the three resulting $(-3)$ curves
and call the result $\overline{IV}$.
This has three singularities of type $\frac{1}{3}(1,1)$.

\item $II^*=\tE_8, III^*=\tE_7$ and $IV^*=\tE_6$:
contract all curves except the central one
and call the result $\overline{R}^*$ for $R=II,III,IV$.
This has three singularities, all of type $A$.
\end{enumerate}

Now suppose that $Y_i\to C_i$ are smooth simply connected
Jacobian elliptic surfaces, each of which has a fibre
$\phi_i$ (with its structure as a scheme)
of type $\tD_4$ over $i$. Let $\pi_i:Y_i\to \bY_i$
denote the contraction of the four $(-2)$-curves
at the extremities, so that $\bY_i\to C_i$
has a fibre $\bphi_i$, with its structure as a scheme, of type $\bD_4$.
So $\bphi_i$ is a copy of $\P^1$ but with multiplicity $2$. 
Say $p_g(Y_i)=h_i$.

\begin{lemma} For each $i$ there is a smooth
$\Z/2$-curve $\tC_i$
and a smooth proper Deligne--Mumford surface
$\tY_i$ with a projective 
(in particular, representable) morphism $\tY_i\to\tC_i$
such that 
\begin{enumerate}
\item the geometric quotients $[\tC_i]$
and $[\tY_i]$ are
$[\tC_i]=C_i$ and $[\tY_i]=\bY_i$,

\item the quotient morphism $\tC_i\to C_i$ is an isomorphism
outside $i$ and the stabilizer group at the unique point $\tii$ of $\tC_i$
that lies over $i$
is $\Z/2$,

\item $\tY_i\to\tC_i$ is smooth over $\tii$ and the fibre $\tphi_{\tii}$
  over $\tii$
is $E/(-1)$ for some elliptic curve $E$,

\item the quotient morphism
$\rho_i:\tY_i\to\bY_i$ is an isomorphism outside the four points
that map to the four nodes on $\bY_i$, where the stabilizer group
is $\Z/2$,

\item the morphism $\tY_i\to\bY_i\times_{C_i}\tC_i$
is finite and birational and

\item the induced morphisms $\tC_i\to\sEll$
are ramified at the points $\tii$ and so
are isomorphic to first order at $\tilde{a},\tilde{b}$.
\end{enumerate}
\begin{proof} 
Locally, each $A_1$ singularity on $\bY_i$ is defined
by an equation $z^2=ty$, where $t$ is a local co-ordinate on $C_i$.
So, if $B_i\to C_i$ is a double cover that
is branched at $i$, then $B_i$ has a local co-ordinate $s$
such that $t=s^2$ and then
$\bY_i\times_{C_i}B_i$ is defined by the equation
$z^2=s^2y$ and the normalization of this is smooth.
The group $\Z/2$ acts on $B_i$,
and then we construct $\tC_i$ locally 
in such a way that near $i$ it is $B_i$ while 
away from $i$ it is $C_i$.
Then $\tC_i=B_i/(\Z/2)$ and $\tY_i$ is the normalization
of $\bY_i\times_{C_i}\tC_i$.
\end{proof}
\end{lemma}

\begin{lemma}\label{identify}
\part[i] There are natural isomorphisms
\begin{eqnarray*}
\pi_i^*\omega_{\bY_i}&\buildrel\cong\over\to&\Omega^2_{Y_i},\\
\pi_i^*\omega_{\bY_i}(\bphi_i)&\buildrel\cong\over\to&\Omega^2_{Y_i}(\phi_i),\\
\rho_i^*\omega_{\bY_i}&\buildrel\cong\over\to&\Omega^2_{\tY_i}\ \textrm{and}\\
\rho_i^*\omega_{\bY_i}(\phi_i)&\buildrel\cong\over\to&\Omega^2_{\tY_i}(2\tphi_{\tii}).
\end{eqnarray*}

\part[ii] These isomorphisms induce isomorphisms
\begin{eqnarray*}
\pi_i^*:H^0(\bY_i,\omega_{\bY_i})&\buildrel\cong\over\to& H^0(Y_i,\Omega^2_{Y_i}),\\
\pi_i^*:H^0(\bY_i,\omega_{\bY_i}(\bphi_i))&\buildrel\cong\over\to& H^0(Y_i,\Omega^2_{Y_i}(\phi_i)),\\
\rho_i^*:H^0(\bY_i,\omega_{\bY_i})&\buildrel\cong\over\to& H^0(\tY_i,\Omega^2_{\tY_i})
\ \textrm{and}\\
\rho_i^*:H^0(\bY_i,\omega_{\bY_i}(\bphi_i))&\buildrel\cong\over\to& H^0(\tY_i,\Omega^2_{\tY_i}(2\tphi_{\tii})).
\end{eqnarray*}
\begin{proof}
This is a straightforward (and well known) 
local calculation for $A_1$ singularities
that depends on two facts: one is that
the non-trivial element $\iota$ of a stabilizer group at a point on
$\tY_i$ that lies over a singularity on $\bY_i$ acts trivially
on a local generator $d\tz_i\wedge dw_i$ of $\Omega^2_{\tY_i}$
where $\tz_i$ is a local co-ordinate on $\tC_i$
such that $\iota^*\tz_i=-\tz_i$ and $w_i$ is a fibre co-ordinate,
and the other is that $Y_i\to\bY_i$ is a crepant resolution.
\end{proof}
\end{lemma}
Now assume that the two $\bD_4$ fibres are
isomorphic. That is, two cross-ratios are equal.
\begin{corollary}\label{8.12} Choose local co-ordinates
$\tz_i$ on $\tC_i$ at $\tii$ such that
$\iota^*\tz_i=-\tz_i$
where $\iota$ is a generator of the stabilizer 
of the point $\tii$. 

\part[i] In terms of these local co-ordinates
there is a Fay plumbing of the stacks $\tY_i$ modulo $t^2$
via the morphisms
$\tC_i\to\sEll$
to give $\tsY'\to\Delta'$.
\part[ii] If $\tsY\to\Delta$ is any lifting
of $\tsY'\to\Delta'$ then $\tsY$
is a smooth $3$-dimensional Deligne--Mumford
stack, whose only non-trivial
stabilizers are copies of $\Z/2$
at each of the four $2$-torsion points
of $\tY_a\cap\tY_b$
(which is isomorphic to the quotient $E/\iota$
of an elliptic curve $E$ by $(-1)$).
\part[iii] Let $\bsY=[\tsY]$. Then $\bsY_0=\bY_a\cup\bY_b$
and $\bY_a\cap\bY_b$ is a fibre of type $\bD_4$
on each of them.
\part[iv]\label{no monodromy} $p_g(\tsY_t)=p_g(Y_a)+p_g(Y_b)$
and there is no monodromy on $H^2(\tsY_t,\Z)$.
\begin{proof}
\DHrefpart{i}: note that, since the $\bD_4$ fibres are isomorphic,
the classifying morphisms $\tC_i\to\sEll$ are then isomorphic
to first order, and the plumbing is constructed
just as in Corollary \ref{7.2}. 

\DHrefpart{ii}: the proof 
of smoothness is the same
as in that of Corollary \ref{construct}
and the rest is immediate.

\DHrefpart{iii}: this is clear.

\DHrefpart{iv}: A  holomorphic $2$-form $\omega_t$ on $\tsY_t$
will specialize to a pair $(\omega_a,\omega_b)$ where
$\omega_i$ lies in $H^0(\tY_i,\Omega^2_{\tY_i}(E/\iota))$
Since $\iota$ acts as $(-1)$ on $E$ the residue of $\omega_i$
is zero, so $\omega_i$ is holomorphic.
(In contrast, we had
$p_g(\sX_t)=g_a+g_b+1$ because there was no
constraint imposed by anti-invariance
to force the vanishing of the residue of the
relevant $2$-form.)

The absence of
monodromy is a well known consequence.
\end{proof}
\end{corollary}
Now suppose that there is a lifting $\tsY\to\Delta$, and fix one.
Then there is a $C^0$ collapsing map
$\g:\tsY_t\to\tsY_0$,
in the context
of orbifolds,
unique up to homotopy. There is a totally isotropic sublattice
$L$ of $H_2(\tsY_t,\Z)$ with a basis
$(A_1,...,A_{h_a},A_{h_a+1},...,A_{h_a+h_b})$
such that the image of 
$(A_1,...,A_{h_a})$ under $\g_*$ is a basis of an isotropic lattice
$L_a$ in $H_2(\tY_a,\Z)$ and the image of
$(A_{h_a+1},...,A_{h_a+h_b})$
is a basis of an isotropic lattice
$L_b$ in $H_2(\tY_b,\Z)$.
We shall assume that each $\tY_i$ is in $L_i$-general position
(see Definition \ref{L general}).
\begin{lemma}\label{gen pos survives}
$\tsY_t$ is in $L$-general position for $t\ne 0$.
\begin{proof} If not, then there is a holomorphic $2$-form
$\omega_t$ on $\tsY_t$ such that $\int_A\omega_t=0$
  for every $A\in L$. We can
assume that $\omega_t$ is not divisible by $t$.
By the argument used in the proof of Lemma \ref{no monodromy}
$\omega_t$ specializes to
a pair $(\omega_a,\omega_b)$, where $\omega_i$
is holomorphic on $\tY_i$ and $\int_{A}\omega_i=0$
for every $A\in L_i$. Then $\omega_a$ and
$\omega_b$ both vanish; this means that $\omega_t$
is divisible by $t$, contradiction.
\end{proof}
\end{lemma}
Poincar{\'e} duality on $\tY_i$ identifies 
$H_2(\tY_i,\Z)$ with a subgroup of $H^2(\tY_i,\Z)$
whose quotient is $2$-elementary.
Correspondingly, by Lemma \ref{gen pos survives},
there is a unique normalized basis
$(\omega^{(1)}(t),...,\omega^{(h)}(t))$
of $H^0(\tsY_t,\Omega^2_{\tsY_t})$
such that $\int_{A_l}\omega^{(j)}(t)=\delta^j_l.$

There are holomorphic $3$-forms $\Omega^{(j)}$
on the smooth 3-fold stack $\tsY$ such that
$$\omega^{(j)}(\lambda)=\Res_{\tsY_\lambda}\Omega^{(j)}/(t-\lambda).$$
We can expand $\Omega^{(j)}$ locally, on the inverse image
in $\tsY$ of the plumbing fixture $\tF$, as
$$\Omega^{(j)}=\sum_{m,n\ge 0}c_{m,n}^{(j)}q^mv^ndq\wedge dv\wedge dw,$$
so that, since $t\equiv q^2-v^2\pmod{t^2}$,
$$\omega^{(j)}(t)\equiv -{\half}\sum_{m,n\ge 0}c_{m,n}^{(j)}q^mv^{n-1}dq\wedge dw.$$
(Here and from now on all congruences are taken modulo $t^2$.)

In a \nbd of $\tY_a$ we have $v=q(1-tq^{-2})^{1/2},$
so that 
$$v^{n-1}=q^{n-1}(1-tq^{-2})^{(n-1)/2}\equiv q^{n-1}(1-(n-1)tq^{-2}/2).$$
Therefore $\omega^{(j)}(t)\vert_{\tY_a}\equiv \omega^{(j)}_{\tY_a}+t\eta^{(j)}_{\tY_a}$
where
\begin{eqnarray*}
\omega^{(j)}_{\tY_a}&=&-{\half}\sum c^{(j)}_{m,n}q^{m+n-1}dq\wedge dw\ \textrm{and}\\
\eta^{(j)}_{\tY_a}&=& {\quarter}\sum c^{(j)}_{m,n}(n-1)q^{m+n-3}dq\wedge dw
\end{eqnarray*}
are $2$-forms (the first holomorphic, the second meromorphic)
on $\tY_a$.
On $\tY_a$ we have $q=\tz_a$ and $w=w_a$, so that
\begin{eqnarray*}
\omega^{(j)}_{\tY_a}&=&-{\half}\sum_{p\ge 0}\left(\sum_{m+n=p} c^{(j)}_{m,n}\right)
\tz_a^{p-1}d\tz_a\wedge dw_a\ \textrm{and}\\
\eta^{(j)}_{\tY_a}&=& {\quarter}\sum_{p\ge 0}\left(\sum_{m+n=p}(n-1) c^{(j)}_{m,n}\right)
\tz_a^{p-3}d\tz_a\wedge dw_a.
\end{eqnarray*}
In a \nbd of $\tY_b$ we have $v=-q(1-tq^{-2})^{1/2}$ 
and on $\tY_b$ we have $q=\tz_b$ and $w=w_b$, so that
$\omega^{(j)}(t)\vert_{\tY_b}\equiv \omega^{(j)}_{\tY_b}+t\eta^{(j)}_{\tY_b}$
where 
\begin{eqnarray*}
\omega^{(j)}_{\tY_b}&=&-{\half}\sum_{p\ge 0}\left(\sum_{m+n=p} (-1)^{n-1}c^{(j)}_{m,n}\right)
\tz_b^{p-1}d\tz_b\wedge dw_b\ \textrm{and}\\
\eta^{(j)}_{\tY_b}&=& {\quarter}\sum_{p\ge 0}\left(\sum_{m+n=p}(-1)^{n-1}(n-1) c^{(j)}_{m,n}\right)
\tz_b^{p-3}d\tz_b\wedge dw_b.
\end{eqnarray*}
The involution $\iota$ acts via $\iota^*\tz_i=-\tz_i$ and $\iota^*w_i=-w_i$,
so that, since each $\omega^{(j)}_{\tY_i}$ and $\eta^{(j)}_{\tY_i}$ is $\iota$-invariant,
we need only consider odd values of the index $p$ in the expansions above.
That is,
\begin{eqnarray}\label{stack expansion}
\omega^{(j)}_{\tY_a}&=&-{\half}\sum_{r\ge 0}\left(\sum_{m+n=2r+1} c^{(j)}_{m,n}\right)
\tz_a^{2r}d\tz_a\wedge dw_a,\\
\eta^{(j)}_{\tY_a}&=& {\quarter}\sum_{r\ge 0}\left(\sum_{m+n=2r+1}(n-1) c^{(j)}_{m,n}\right)
\tz_a^{2r-2}d\tz_a\wedge dw_a,\\
\omega^{(j)}_{\tY_b}&=&-{\half}\sum_{r\ge 0}\left(\sum_{m+n=2r+1} (-1)^{n-1}c^{(j)}_{m,n}\right)
\tz_b^{2r}d\tz_b\wedge dw_b\ \textrm{and}\\
\eta^{(j)}_{\tY_b}&=& {\quarter}\sum_{r\ge 0}\left(\sum_{m+n=2r+1}(-1)^{n-1}(n-1) c^{(j)}_{m,n}\right)
\tz_b^{2r-2}d\tz_b\wedge dw_b
\end{eqnarray}
\noindent and $c^{(j)}_{m,n}=0$ if $m+n$ is even. In particular,
these formulae show that each $\eta^{(j)}_{\tY_i}$ has
only a double pole along $\tphi_{\tii}$, and is otherwise holomorphic.

The next lemma is well known. 
It holds for any one-parameter degeneration of surfaces,
not only for the families that we have constructed via Fay plumbings.
However, the results such as Lemma \ref{unique} below that 
follow for the derivatives $\eta^{(j)}_{X_i}$
do not hold in such generality, because for more general 
glueing data a formula such as $t\equiv q^2-v^2{\pmod {t^2}}$
will not hold. It follows that
for a general one-parameter degeneration
it is not possible to control the orders of the poles
of the forms $\eta^{(j)}_{\tY_i}$ (which are the restrictions to $\tY_i$
of the derivatives at $t=0$
of the forms $\omega^{(j)}(t)$) along $\tphi_{\tii}$.

\begin{lemma} 
\part[i] $\omega^{(j)}_{\tY_i}$ is holomorphic on $\tY_i$
for every $j\in [1,h_a+h_b]$.
\part[ii] $\omega^{(j)}_{\tY_b}=0$ for every $j\in [1,h_a]$
and $\omega^{(j)}_{\tY_a}=0$ for every $j\in[h_a+1,h_a+h_b]$.
\part[iii] $(\omega^{(1)}_{\tY_a},...,\omega^{(h_a)}_{\tY_a})$
and $(\omega^{(h_a+1)}_{\tY_b},...,\omega^{(h_a+h_b)}_{\tY_b})$
are bases of the vector spaces $H^0(\tY_a,\omega_{\tY_a})$
and $H^0(\tY_b,\omega_{\tY_b})$, respectively.
\part[iv] These bases are normalized with respect to the given $A$-cycles
on $\tY_a$ and $\tY_b$.
\begin{proof}
This is a restatement of Lemmas \ref{no monodromy} and \ref{gen pos survives}.
\end{proof}
\end{lemma} 

\begin{lemma}\label{unique}\label{up to}
\part[i] Up to scalars there is
a unique meromorphic $2$-form 
$\teta_{\tY_i}\in H^0(\tY_i,\Omega^2_{\tY_i}(2\tphi_{\tii}))$
such that every $\int_{A_l}\teta_{\tY_i}=0$.

\part[ii] Every $\eta^{(j)}_{\tY_i}$ is a multiple of $\teta_{\tY_i}$.
\begin{proof} 
\DHrefpart{i} follows from the facts that
$\dim H^0(\tY_i,\Omega^2_{\tY_i}(2\tphi_{\tii}))=p_g(\tY_i)+1$
and that $\tY_i$ is in $L_i$-general position.
\DHrefpart{ii} is then an immediate consequence
of \DHrefpart{i} and the observation above that
each $\eta^{(j)}_{\tY_i}$ lies in $H^0(\tY_i,\Omega^2_{\tY_i}(2\tphi_{\tii}))$.
\end{proof}
\end{lemma}

We shall refer to the forms $\teta_{\tY_i}$ as being \emph{normalized}
by the requirements that
$\teta_{\tY_i}\in H^0(\tY_i,\Omega^2_{\tY_i}(2\tphi_{\tii}))$
and every $\int_{A_l}\teta_{\tY_i}=0$.

Construct four row vectors
\begin{eqnarray*}
{\underline{\eta}}_{\tY_i}&=&[\eta^{(1)}_{\tY_i},...,\eta^{(h_a+h_b)}_{\tY_i}],\\
{\underline{\omega}}_{\tY_a}&=&[\omega^{(1)}_{\tY_a},...,\omega^{(h_a)}_{\tY_a}]\ \textrm{and}\\
{\underline{\omega}}_{\tY_b}&=&[\omega^{(h_a+1)}_{\tY_b},...,\omega^{(h_a+h_b)}_{\tY_b}],
\end{eqnarray*}
each consisting of meromorphic $2$-forms on the surface indicated.
As an abbreviation, define two further row vectors, each consisting of
numbers, by
$${\underline{\omega}}_{\tY_i}(i)=
\frac{{\underline{\omega}}_{\tY_i}}{d\tz_i\wedge dw_i}\left(i\right).$$

\begin{proposition}\label{important} We can scale the forms $\teta_{\tY_i}$ such that
in terms of the local co-ordinates $\tz_i$ they are given by
\begin{eqnarray*}
\teta_{\tY_a}&=&{\quarter}(\tz_a^{-2}+\ \textrm{h.o.t.})d\tz_a\wedge dw_a\ \textrm{and}\\
\teta_{\tY_b}&=&-{\quarter}(\tz_b^{-2}+\ \textrm{h.o.t.})d\tz_b\wedge dw_b
\end{eqnarray*}
\noindent and
there is an equality of row vectors
$${\underline{\eta}}_{\tY_i}=\teta_{\tY_i}
[{\underline{\omega}}_{\tY_a}(a),-{\underline{\omega}}_{\tY_b}(b)].$$
\begin{proof} 
Examining the first few terms in the expansion
provided by the formulae (\ref{stack expansion}) \emph{et\ seq.} gives
\begin{eqnarray*}
\omega_{\tY_a}^{(j)}(a)&=&{-\half}(c^{(j)}_{1,0}+c^{(j)}_{0,1}),\\ 
\omega_{\tY_b}^{(j)}(b)&=&{-\half}(-c^{(j)}_{1,0}+c^{(j)}_{0,1}),\\
\eta^{(j)}_{\tY_a}&=&{\quarter}(-c^{(j)}_{1,0}\tz_a^{-2}+\textrm{h.o.t.})d\tz_a\wedge dw_a\ \textrm{and}\\
\ \eta^{(j)}_{\tY_b}&=&{\quarter}(c^{(j)}_{1,0}\tz_b^{-2}+\textrm{h.o.t.})d\tz_b\wedge dw_b
\end{eqnarray*}
for all $j=1,...,h_a+h_b$.

If $j\le h_a$ then $\omega^{(j)}_{\tY_b}=0$, so that
$c^{(j)}_{1,0}=c^{(j)}_{0,1}=-\omega^{(j)}_{\tY_a}(a).$

If $j\ge h_a+1$ then
$\omega^{(j)}_{\tY_a}=0$, so that
$c^{(j)}_{1,0}+c^{(j)}_{0,1}=0$
and $\omega^{(j)}_{\tY_b}(b)=c^{(j)}_{1,0}$.
The result is proved.
\end{proof}
\end{proposition}
Via the identifications of Lemma \ref{identify} we also regard these as vectors of forms
on $Y_i$ and on $\bY_i$ and also write them as ${\underline{\eta}}_{Y_i}$, etc.

Let $\Psi_{Y_i}$ denote the period matrix of $Y_i$,
normalized
with respect to the $2$-cycles $A_1,...,A_{h_i}$
above. Recall that normalizing means that
$\int_{A_j}\omega^{(k)}_{Y_i}=\delta_j^k$ for
each $i$ and for
the appropriate values of $j,k$
and that $\Psi_{Y_i}$ is an $h_i\times (11h_i+8)$ matrix
while $\Psi_t:=\Psi_{\tsY_t}$ is an $h\times (11h+8)$ matrix.
Note, however, that each $Y_i$ contains four disjoint $(-2)$-curves
(components of the $\tD_4$-fibre that is designated on each),
so that each $\Psi_{Y_i}$ contains an $h_i\times 4$ block of zeroes,
which can be discarded to give matrices
that we shall continue to denote by $\Psi_{Y_i}$,
each of whose shape is $h_i\times(11h_i\times 4)$.
\begin{proposition}\label{8.2}\label{11.8}
\part[i] 
The normalized period matrix $\Psi_t$ of $\tsY_t$ is given,
after a suitable re-arrangement of the blocks that comprise it, by
$$\Psi_t=
  \left[
  \begin{array}{cc}
    {\Psi_{Y_a}} & 0\\
    0 &{\Psi_{Y_b}}
  \end{array}
  \right]
+
t[\underline{\omega}_{Y_a}(a),-\underline{\omega}_{Y_b}(b)]\otimes
[\underline{I}_a,\underline{I}_b]+{\textrm{h.o.t.}}$$
where ${\underline{I}}_i$ is a vector of integrals of 
$\teta_{Y_i}$ 
over cycles on $Y_i$.

\part[ii]
The derivative $(d\Psi_t/dt)\vert_{t=0}$
of the period matrix of $\sY_t$ at $t=0$ is the rank $1$ matrix
$[\underline{\omega}_{Y_a}(a),-\underline{\omega}_{Y_b}(b)]\otimes
[\underline{I}_a,\underline{I}_b].$
\begin{proof} This is an immediate consequence
of Proposition \ref{important} and the fact that
$$\omega^{(j)}(t)\vert_{\tY_i}\equiv\omega^{(j)}_{\tY_i}+t\eta^{(j)}_{\tY_i},$$
modulo $t^2$, for each $j=1,...,h_a+h_b$ and for each $i=a,b$.

Note that the vector $[\underline{I}_a,\underline{I}_b]$
can be seen to have the correct length
(so that the matrix that describes $(d\Psi_t/dt)\vert_{t=0}$
has the correct shape) by omitting the zeroes, four in each
$[\underline{I}_i]$, that arise from integrating
$\teta_i$ around the four exceptional $(-2)$ curves
on $Y_i$ that are contracted in $\bY_i$.
\end{proof}
\end{proposition}

\begin{lemma}\label{non-zero}
Suppose that $h_i\ge 1$ and that 
$Y_i$ is \emph{either} a special elliptic surface \emph{or}
generic. Then the class $[\teta_{Y_i}]$ in $H^2(Y_i,\C)$
is non-zero.
\begin{proof}
Suppose first that $Y_i$ is special. Then $Y_i$ is birational to 
the geometric quotient $C\times E/\iota$, where $C$ is hyperelliptic
of genus $h_i$. Suppose that $P\in C$ is a fixed point of $\iota\vert_C$,
so that the image of $\{P\}\times E$ in $Y_i$ is a $\tD_4$ fibre $\phi_i$.

There is an $\iota$--anti-invariant section $\s$ of $H^0(C,\Omega^1_C(2P))$
that is not holomorphic; modulo the subspace $H^0(C,\Omega^1_C)$
this section is unique up to scalars. Moroever, $\s$ has no residues
and defines a non-trivial
class in $H^{0,1}(C)$. 
Choose a non-zero section $\tau$ of $H^0(E,\Omega^1_E)$.
Then $\tau$ is $\iota$--anti-invariant, and $\s\wedge\tau$
defines a non-zero section of $H^0(Y_i,\omega_{Y_i}(\phi_i))/H^0(Y_i,\omega_{Y_i})$.
Inspection of the leading term of $\teta_{Y_i}$ shows
that $\teta_{Y_i}$ defines a non-zero section of the same
$1$-dimensional vector space, so that $\s\wedge\tau$ is equivalent to
$\lambda\teta_{Y_i}$ modulo $H^0(Y_i,\omega_{Y_i})$, for some $\lambda\in\C$.
But $\s\wedge\tau$ defines a non-zero class in $H^{1,1}(Y_i)$,
and the lemma is proved when $Y_i$ is special.

The result for generic $Y_i$ follows at once.
\end{proof}
\end{lemma}

We can now re-state Proposition \ref{8.2} in intrinsic terms.
We shall make use of Remark \ref{transv} below,
to the effect that $[\teta_{Y_i}]$ in fact lies in $Fil^1(H^2(Y_i,\C))$. 

\begin{corollary}\label{rank 1} $(d\Psi_t/dt)\vert_{t=0}$
is a linear map
$$H^0(Y_a,\omega_{Y_a})\oplus H^0(Y_b,\omega_{Y_b})\to H^{1,1}(Y_a)
\oplus H^{1,1}(Y_b)$$
whose rank is at most $1$.

If at least one of the surfaces $Y_a$ or $Y_b$ 
is special or generic, then the rank of $(d\Psi_t/dt)\vert_{t=0}$ equals $1$, 
its kernel is the hyperplane
$\{(\s_a,\s_b)\vert \s_a(a)+\s_b(b)=0\}$ 
in $H^0(Y_a,\omega_{Y_a})\oplus H^0(Y_b,\omega_{Y_b})$
and its image is the line
spanned by the element $([\teta_{Y_a}],[\teta_{Y_b}])$.
\noproof
\end{corollary}

\begin{remark}
\part[i]\label{transv}
Recall also that the final $h\times h$ block (with respect to an appropriate
choice of basis) of the period matrix,
and so of the derivative $(d\Psi/dt)\vert_{t=0}$, is skew-symmetric. But a skew-symmetric
matrix of rank $1$ vanishes identically; for $(d\Psi/dt)\vert_{t=0}$ 
this is predicted by Griffiths transversality.
In terms of the bases that we have used, this means that the final piece of
length $h_i$ in
each vector $\underline{I}_i$ vanishes. Stated in intrinsic terms,
this means that each class $[\teta_{Y_i}]$ lies in $H^{1,1}(Y_i)$.

\part[ii]
Suppose that $Y_b$ is rational.
Then $(d\Psi_t/dt)\vert_{t=0}$ is given by
$$(d\Psi_t/dt)\vert_{t=0}=[\underline{\omega}_{Y_a}(a)]\otimes 
[\underline{I}_a,\underline{I}_b].$$ 
Note that $Y_b$ has $4$ moduli if the constraint
that its $\tD_4$ fibre $\phi_b$ should be isomorphic to $\phi_a$
is ignored. This constraint reduces the number of moduli of $Y_b$ to $3$.

There is a desingularization $\sY\to[\tsY]$
of the geometric quotient 
such that the closed fibre $\sY_0$ of $\sY\to\Delta$
is semi-stable and
has six components: two are the surfaces $Y_i$ for $i=a,b$
and four are
copies of $\P^2$ each of whose normal bundle is
$\sO(-2)$.
So $\sY\to \Delta$ has a birational model $\sY'\to\Delta$ with good
reduction, which is obtained from $\sY$ after flopping four times
(each time in a curve in $Y_b$ that is a section of 
the elliptic fibration and meets one of the $V_i$),
then contracting the strict transforms of the $V_i$
to curves
and finally contracting the strict transform of $Y_b$
to a curve.
The closed fibre $\sY'_0$ is isomorphic to $Y_a$.
That is, a Fay plumbing gives a one-parameter
deformation of a surface $Y_a$
where the derivative of the period matrix is given
explicitly, provided that $Y_a$ contains a 
$\tD_4$-fibre. Since $Y_b$ has three moduli
the vector ${\underline{I}}_b$ has three moduli,
so $(d\Psi_t/dt)\vert_{t=0}$ has three moduli.

This is an analogue of Example 3.5
on p. 45 of \cite{F1}.
\end{remark}
The next lemma is needed for the proof of Proposition \ref{12.9}.

\begin{lemma}\label{root}
  If $\delta$ is a vertical $(-2)$-curve on $Y_i$
  that does not lie in the designated $\tD_4$-fibre $\phi_i$,
  then $\int_\delta\teta_{Y_i}=0$.
  \begin{proof}
    We can take $i=a$.

    By construction, the first order plumbing $\sY'\to\Delta'$ is trivial
    along $Y_a-T_a$, where $T_a$ is a suitable \nbd of $\phi_a$.
    So $\int_\delta\omega^{(j)}(t)\equiv 0{\pmod{t^2}}$
    for $j=1,...,h_a$.
    But $\omega^{(j)}(t)\vert_{Y_a}= \omega^{(j)}_{Y_a}+t\eta^{(j)}_{Y_a}+{h.o.t.}$
    and the lemma follows.
  \end{proof}
\end{lemma}
  
Now suppose that 
$f_1:Y_1\to C_1,...,f_r:Y_r\to C_r$ are smooth simply connected Jacobian 
elliptic surfaces such that $Y_i\to C_i$ has
$\tD_4$-fibres over a finite non-empty 
set $\{P_{ij}\}$ of points in $C_i$. Let $Y_i\to \bY_i$
denote the contraction of each of these designated
$\tD_4$-fibres to a $\bD_4$-fibre. 

Suppose also that
$\G$ is a connected tree with $r$ vertices
and that we can associate the surfaces $\bY_i$ to the vertices
of $\G$ such that two vertices $i,j$ are joined in $\G$
if and only if the $\bD_4$-fibre on $\bY_i$ that lies over $P_{ij}$
is isomorphic to the $\bD_4$-fibre on $\bY_j$ that lies over $P_{ji}$.

Then the plumbings modulo $t^2$ of the stacks $\tY_i$ 
that has
just been described can be iterated to give
$$\sY'=\sY'_\G\to\Delta'_{r-1}=\Sp\C[t_1,...,t_{r-1}]/(t_1^2,...,t_{r-1}^2)$$
where the closed fibre is $\sum\tY_i$ arranged
according to the tree $\G$.
Let $\underline{\omega}_{Y_i}$
denote a normalized basis of $H^0(Y_i,\Omega^2_{Y_i})$
and $\teta_{ij}$ a normalized element of
$H^0(Y_i,\Omega^2_{Y_i}(\phi_{ij}))$.

\begin{proposition}\label{12.9} 
Assume that $\sY'\to \Delta'_{r-1}$
can be lifted to a family
$\sY\to S_{t_e}$ over an $(r-1)$-dimensional
polydisc. 

\part[i] Then there is no monodromy on $H^2$ of the geometric generic fibre
$\sY_t$ and the period
matrix $\Psi(\sY_t)$ is given by
$$\Psi(\sY_t)=[\Psi(Y_1),...,\Psi(Y_r)]+
\sum_e t_e\Pi_e +{h.o.t.},$$
where $\Pi_e$ is the rank $1$ matrix
$$[{\underline{\omega}}_{Y_i}(P_{ij}),
-{\underline{\omega}}_{Y_j}(P_{ji})]\otimes 
[{\underline{I}}_{ij},{\underline{I}}_{ji}]$$
and ${\underline{I}}_{ij}$ is a vector of integrals of $\teta_{ij}$ around $2$-cycles on
$Y_i$.

\part[ii] If each $Y_i$ is either special or generic, then the image of the
tangent space $T_0S_{t_e}$ under the derivative of the period map
is of dimension $r-1$ and is the vector space spanned by the matrices
$\Pi_e$.
\begin{proof} \DHrefpart{i} is a straightforward consequence
  of Proposition \ref{8.2}. Lemma \ref{root} tells us that,
  when considering integrals of $\teta_{ij}$ around cycles
  on $Y_i$, we can omit the $(-2)$ curves that lie in
  the $\tD_4$-fibres on $Y_i$ and are disjoint from $\phi_{ij}$.
  This makes it possible to write the vectors ${\underline{I}}_{ij}$
  in such a way that the matrices $\Pi_e$, when suitably enlarged
  by adding various blocks of zero matrices, lie
  in the same vector space of matrices.

  \DHrefpart{ii} then follows from Lemma \ref{non-zero}.
\end{proof}
\end{proposition}
We shall see, in Lemma \ref{sm is smooth}, that
the liftings $\sY\to S_{t_e}$ exist when each surface $Y_i$ is
a special elliptic surface. Then,
in Theorem \ref{8.y},
which will lead to the main result, Theorem \ref{main},
we shall take all the surfaces $Y_i$
to be K3 surfaces and $\G$ will be an alkane.
\end{section} 
\begin{section}{Stable reduction of surfaces}
Chakiris proved \cite{C1}, \cite{C2} a generic Torelli theorem for
Jacobian elliptic surfaces over $\P^1$ by reducing the problem
to special elliptic surfaces.

To achieve this reduction he extended the domain of the
period map to make the map proper.
In turn, he did this by first showing
that a one-parameter degeneration of Jacobian elliptic surfaces over $\P^1$
without monodromy can be put into a certain standard form.
See, e.g., the statement $(**)$ in \cite{C2}, top of p. 174
or the ``stable reduction'' theorem on p. 231 of \cite{C1}.
Note, however, that this version of a stable reduction theorem
gives a closed fibre that contains curves of cusps,
and so is not semi log canonical (slc) in the sense of the
MMP. 

In this section 
we shall refine his result, so that the period map 
becomes proper over each of the loci $\sW_{h_1,...,h_r}$
defined in Definition \ref{defn of sW}.
For this, we use his ideas strengthened by the MMP, which was not
available to him.  
The main result is Theorem \ref{10.x}.

We say that a special elliptic surface $X$ is \emph{sesqui-special}
if it is associated to a product $E\times C$ where the
Jacobians $E=\Jac(E),\ \Jac(C)$ of $E$ and $C$ have no multiplication
and $\Hom(E,\Jac(C))=0$. If $C$ and $E$ are generic
then $X$ is sesqui-special \cite{Z}. A very general point
in $\sW_{h_1,...,h_r}$ is defined by a configuration
of sesqui-special surfaces $X_1,...,X_r$ where $X_i$
is associated to $E\times C_i$ and also
$\Hom(\Jac(C_i),\Jac(C_j))=\Z\delta_{ij}$.

\begin{theorem}\label{10.x} Suppose that $\sX\to\Delta$ is a
$1$-parameter
degeneration of simply connected Jacobian elliptic surfaces
of geometric genus $h\ge 1$
which is semi-stable (in the usual sense that the closed fibre
$\sX_0$ is reduced with normal crossings)
and that there is no monodromy on the cohomology
of the geometric generic fibre $\sX_{\bar\eta}$.
Assume also that, under the period map,
the image of $0\in\Delta$ is the direct sum
of the period matrices of sesqui-special
elliptic surfaces
$\tV_1,...,\tV_r$ that define a point in
$\sW_{h_1,...,h_r}$.

Then there is a birationally equivalent
model $\sY\to\Delta$ with the following properties:

\part[i] $\sY$ has 
$\Q$-factorial canonical singularities;

\part[ii] the closed fibre $\sY_0$ has slc singularities;

\part[iii] the irreducible components 
of $\sY_0$ are the singular models $V_i$ of the $\tV_i$
with $\bD_4$-fibres;

\part[iv] if $V_i\cap V_j$
is not empty then it is a copy of $\P^1$
and contains $4$ points at all of which $V_i$ and $V_j$ each has a node;

\part[v] each triple intersection $V_i\cap V_j\cap V_k$ is empty;

\part[vi] $\sY_0$ is formed by arranging the surfaces $V_i$ in a tree.
\begin{proof}
  To begin, suppose that $\sX\to\Delta$ is an arbitrary
  semi-stable family
of smooth minimal elliptic surfaces of Kodaira dimension $1$.

Then run the MMP in two steps, as follows.
\begin{enumerate}
  \item
        Run a $K_{\sX/\Delta}$ MMP on $\sX\to\Delta$ and let 
$\sX_1\to\Delta$ be the output. Then $\sX_1$ has $\Q$-factorial
terminal singularities and $K_{\sX_1/\Delta}$ is semi-ample
So some relative pluricanonical system
$\vert mK_{\sX_1/\Delta}\vert$ defines an algebraic fibre space 
$f:\sX_1\to S$ where $S\to\Delta$ is a semi-stable
family of curves (so that $S$ has singularities of type $A$)
and $K_{\sX_1/\Delta}$ pulls back from an ample $\Q$-line bundle on $S$. 
Moreover, the closed fibre $\sX_{1,0}$ has slc
singularities. Replace $\sX$ by $\sX_1$.

\item
If there are surfaces $E_i$ in $\sX$ such that
$f(E_i)$ is a point then there are only finitely many such.
For suitable $\a_i\in\Q$ with $0<\a_i\ll 1$ run a
$(K_{\sX/S},\sum\a_iE_i)$ MMP on $\sX\to S$.
The output is a birational map
$\sX-\to\sX_1$ under which the strict transform of each $E_i$ 
is of dimension at most $1$.
Since $K_{\sX}$ is trivial in a neighbourhood of $\sum E_i$,
the rational map $\sX-\to\sX_1$ is regular outside $\sum E_i$,
so that
the rational map $\sX_1-\to S$ is a morphism,
all its fibres are $1$-dimensional, $\sX_1$ has 
canonical singularities
and $\sX_{1,0}$ has slc singularities. 
Moreover, $\sX_1$ is $\Q$-factorial, by general properties of the MMP.
Replace $\sX$ by $\sX_1$.
\end{enumerate}
At this stage of the argument $\sX$ has $\Q$-factorial canonical singularities
and $\sX_0$ has slc singularities.
Moreover, there is a surface $S$, a semi-stable morphism $g:S\to\Delta$
(so that $S$ has singularities of type $A$) and a morphism $f:\sX\to S$
with only one-dimensional fibres
such that $K_{\sX/\Delta}$ is the pullback under $f$
of a $g$-ample $\Q$-line bundle on $S$.

Now assume that there is no monodromy on the cohomology $H^2(\sX_{\bar\eta})$
of the geometric generic fibre.
This is equivalent to $p_g(\sX_{\bar\eta})=\sum p_g(\tV)$,
where $\tV$ runs over the minimal resolutions of the components
of $\sX_0$. Assume also that $\sX_{\bar\eta}$
is simply connected and that the generic fibre
$\sX_\eta$ is Jacobian. However, we shall make no assumption
about the image of $0$ under the period map
until after Lemma \ref{9.8+2}.

Say $S_0=\sum C_i$, $X_i=f^{-1}(C_i)$,
$f_i:X_i\to C_i$ the restriction of $f:\sX\to S$,
$\nu_i:X_i^\nu\to X_i$ the normalization, $\tX_i\to X_i^\nu$ the
minimal resolution. We have $K_{\sX/\Delta}\vert_{X_i}\sim f_i^*(\a_i)$
for some $\a_i\in\Q$ with $\a_i>0$. Say $C_i\cap C_i=P_{ij}$
if the intersection is not empty.

\begin{lemma} \part[i] Each $C_i$ is isomorphic to $\P^1$.
  \part[ii]\label{has section} $f_i:X_i\to C_i$ has a section.
  \begin{proof} $S_0$ is a specialization of $\P^1$,
    and \DHrefpart{i} follows.

    \DHrefpart{ii}: Suppose that $D\subset\sX$ is the
    Zariski closure of the given generic section.
    Then $D\cap X_i$ consists of a section with also
    some vertical components.
  \end{proof}
\end{lemma}

\begin{lemma}\label{slc class}\label{9.6+1} 
$X_i$ and $X_j$ are Cohen--Macaulay 
  and are smooth at each generic point of $f^{-1}(P_{ij})$.
\begin{proof} This is a consequence of the
  classification of slc singularities. In particular,
    smoothness at each generic point of $f^{-1}(P_{ij})$
  follows from the fact that slc singularities
  have normal crossings in codimension $1$.
\end{proof}
\end{lemma}

By the classification of slc singularities,
the generic fibre of each $f_i$ is elliptic or a rational cycle.
For each $i$, put $C_i^0=C_i-\Sing(S_0)$.

\begin{lemma}\label{9.3} The maps
$H^2(\sX_0,\sO)\to\oplus H^2(\tX_i,\sO)$ and
$H^2(X_i^\nu,\sO)\to H^2(\tX_i,\sO)$ are all isomorphisms.
\begin{proof} $H^2(\sX_0,\sO)\to\oplus H^2(\tX_i,\sO)$
is surjective. Since the formation of the groups $H^i(\sX_t,\sO)$
commutes with specialization the lemma
follows from the assumption that 
$p_g(\sX_{\bar\eta})=\sum p_g(\tX_i)$.
\end{proof}
\end{lemma}

\begin{lemma}\label{9.11}
  Suppose that $X_i=f^{-1}(C_i)\to C_i$ is generically smooth
  and that $h^2(X_i,\sO)>0$.
  Then $X_i$ is normal with rational singularities,
  and over $C_i^0$ it is Gorenstein.
\begin{proof} Let $Q\in C_i^0$
and let superscript $h$ denote henselization at $Q$.
Note that $S_0^h=C_i^h$.
Then $\sX^h-X_i^h\to S^h-S_0^h$
is a Jacobian elliptic fibration, so that
$K_{\sX^h-X_i^h}$ pulls back from a line bundle
on $S^h-S_0^h$ and so is trivial. Since $X_i^h$ is an
irreducible and principal Weil divisor in $\sX^h$,
the restriction homomorphism $\Cl(\sX^h)\to\Cl(\sX^h-X_i^h)$
of Weil divisor class groups
is an isomorphism. Therefore $K_{\sX^h}$ is trivial.
So $\sX$ is Gorenstein along $f^{-1}(Q)$ and $K_{\sX^h}$
pulls back from a line bundle on $S^h$.
Therefore $X_i^h$ is also Gorenstein
and $\omega_{X_i^h}$ also
pulls back from a line bundle on $C_i^h$.

There is a commutative diagram
$$\xymatrix{
  &{\tX_i}\ar[dl]_{k_i}\ar[dr]^{g_i}\\
  {X_i^{min}}\ar[dr]_{f_i^{min}}&&{X_i}\ar[dl]^{f_i}\\
  &{C_i}
}
$$
where
$g_i:\tX_i\to X_i$ is the minimal resolution,
$k_i:\tX_i\to\tX_i^{min}$ is the minimal model relative to $C_i$
and $f_i^{min}:\tX_i^{min}\to C_i$ is the resulting morphism.
Then
$$f_i^{min*}K_{\tX_i^{min}}+E\sim K_{\tX_i}\sim g_i^{*}K_{X_i}-Z$$ for some
non-negative integral divisors $Z,E$.
Moreover, $Z>0$ if $X_i$ is not normal along $f_i^{-1}(Q)$,
since there is a non-zero contribution to $Z$ made by the conductor ideal.
But $k_i^*K_{X_i^{min}}\sim g_i^*K_{X_i}-(Z+E)$
then shows that $h^2(X_i^{min},\sO_{X_i^{min}})< h^2(X_i,\sO_{X_i}),$
contradiction.
So $X_i$ is normal along $f_i^{-1}(Q)$, so that over $C_i^0$
it is
normal and Gorenstein.

That $X_i$ is normal everywhere follows from Lemma \ref{slc class};
it remains to prove that $X_i$ has rational singularities.

So suppose that there is at least one irrational singularity.
We know that $H^2(X_i,\sO_{X_i})\to H^2(\tX_i,\sO_{\tX_i})$
is an isomorphism; then $\tX_i$ is elliptic over $\P^1$
with a section
and $p_g(\tX_i)>0$, so that, by considering the Albanese
variety of $\tX_i$, we see that $H^1(\tX_i,\sO_{\tX_i})=0$.
Via the Leray spectral sequence
$$E^{pq}_2=H^p(X_i,R^qg_{i*}\sO_{\tX_i})\Rightarrow H^{p+q}(\tX_i,\sO_{\tX_i})$$
we then see that $R^qg_{i*}\sO_{\tX_i}=0$, and the singularities are rational.
\end{proof}
\end{lemma}

Suppose that $C,D$ are irreducible components of $S_0$, $C\cap D=P$,
$Y=f^{-1}(C)$, $Z=f^{-1}(D)$. Since $S$ has an isolated
singularity of type
$A$ at $P$, we can write the henselization
$S^h$ of $S$ at $P$ as a geometric quotient
$S^h=[S'/(\Z/N)]$ where $S'$ is smooth and local,
$S'\to\Delta$ is semi-stable and the closed fibre
$S'_0$ has two smooth branches $C',D'$.
Moreover, $\Z/N$ acts with opposite weights
on these two branches, since $S^h$ is of type $A$.
In particular, it acts effectively on each of them.

Let $\sX'$ denote the normalization
of $\sX\times_S S'$, with induced morphism
$f':\sX'\to S'$; this morphism is proper.
Let $P'\in S'$ be the point lying over $P$.

Set $\sX^h=\sX\times_SS^h$ and let $f^h:\sX^h\to S^h$
denote the induced morphism.

Since $\Z/N$ acts freely in codimension $1$ on $S'$, it also
does so on $\sX'$. So the quotient map
$$\sX'\to [\sX'/(\Z/N)]=\sX^h$$
is {\'e}tale outside the $1$-dimensional locus $f^{-1}(P)$.
Since $\sX^h$ has canonical singularities, so does $\sX'$. In particular,
$\sX'$ is Cohen--Macaulay. Since $f':\sX'\to S'$ is dominant with
equi-dimensional fibres and $S'$ is regular,
it follows that $f'$ is flat.

Suppose that $\Delta\to\Delta$ is a finite base change.
Then the same argument shows that
$\sX'\times_\Delta\Delta$ also has canonical singularities,
and so $\sX'_0$ has slc singularities.

The group $\Z/N$ acts effectively on the branches $(C',P')$
and $(D',P')$, and $[C'/(\Z/N)]=C^h$,
$[D'/(\Z/N)]=D^h$.
So $\sX_0'=Y'\cup Z'$, where
$Y'=f'^{-1}(C')$, $Z'=f'^{-1}(D')$.

Let $p:Y'\to C'$ and $q:Z'\to D'$ denote the induced morphisms;
they are proper, and $Y',\ Z'$ are partial
normalizations of $Y\times_CC',\ Z\times_DD'$, respectively.

Since $Y,Z$ are Cohen--Macaulay, we can identify
$Y^h:=Y\times_CC^h=[Y'/(\Z/N)]$
and $Z^h:=Z\times_DD^h=[Z'/(\Z/N)]$.

Since $C',D'$ are principal divisors on $S'$, the divisors $Y'$ and $Z'$ are 
principal on $\sX'$.

Extend the $ADE$ notation in the usual way,
to include $A_0=\A^2$,
$A_\infty=(xy=0)$ and $D_\infty=(x^2=y^2z)$.

Suppose that $\xi$ is a closed point of $f'^{-1}(P')$.
Henselize $\sX',Y',Z'$ at $\xi$ to get a 3-fold germ $\sX''$
and principal divisors $Y'',Z''$ on it.
The stabilizer of $\xi$ is a subgroup $H\cong\Z/M$
of $\Z/N$ that acts freely in codimension $1$ on $\sX''$.
Put $\sX^{loc}=[\sX''/H]$, $Y^{loc}=[Y''/H]$, $Z^{loc}=[Z''/H]$;
these are localizations of $\sX,Y$ and $Z$, respectively.

\begin{lemma}\label{previous}
  \label{9.4}
  \begin{enumerate}
    \item There are just two possibilities:
  \begin{enumerate}
    \item
  \emph{either} $\sX''$ is smooth, $\sX''_0=A_\infty$,
$Y''\cap Z''$ is smooth and
$\sX_0^{loc}=[A_\infty/\frac{1}{M}(1,-1,1)]$
where $M\in\{1,2,3,4,6\}$,
\item \emph{or}\label{(B)}
$\sX''_0$ is a degenerate cusp of multiplicity
at most $4$, $Y'',Z''$ are of type $A$,
$Y''\cap Z''$ is a nodal curve, $M\in\{1,2\}$
and $\sX_0^{loc}$ is either (A) a degenerate cusp of multiplicity ${}\le 4$
or (B) the geometric quotient of such a degenerate cusp by $\Z/2$.
\end{enumerate}

\item In case (\ref{(B)}), $Y^{loc}$ and $Z^{loc}$ are both of type $A_{\ne\infty}$
and $Y^{loc}\cap Z^{loc}$ is a smooth curve.

\item $\sX''$ and $\sX''_0$ are Gorenstein.
  \end{enumerate}
\begin{proof}
We know that $\sX''$ is canonical, $\sX''_0=Y''\cup Z''$ is slc
and that $Y'',Z''$ are principal divisors on $\sX''$.

Next, the classification of slc singularities with at least
two branches shows that the curve $Y''\cap Z''$
is either smooth or a plane node. Then 
in the first case $\sX''$ is smooth,
$\sX''_0=A_\infty$ and 
$\sX_0^{loc}=A_\infty/\frac{1}{M}(a,-a,1)$,
where $a$ is prime to $M$,
while in the second case $\sX''_0$ is a degenerate cusp
and $Y'',Z''$ are of type $A$.

Because $p:Y'\to C'$ is a Jacobian semi-stable family of elliptic curves
and $H$ acts effectively on $Y'$, the classification of automorphism
groups of elliptic curves shows that $M\in\{1,2,3,4,6\}$.
If the fibre $Y'\cap Z'$ over $P'$ of $Y'\to C'$ is singular, then 
$M\in\{1,2\}$. In either case we can take $a=1$.

For \DHrefpart{ii} we can suppose that $M=2$ and $H=\langle\iota\rangle$.

Suppose that $\iota$ switches the two branches of $Y''\cap Z''$.
Then each of $Y^{loc}=[Y''/\iota]$ and
$Z^{loc}=[Z''/\iota]$ is of type $A$ or $D$,
and $Y^{loc}\cap Z^{loc}$ is smooth. The classification
of non-isolated slc singularities then shows
that $Y^{loc}$ and $Z^{loc}$
are both of type $A$.

If $\iota$ fixes the two branches then
$\sX_0^{loc}$ is a degenerate cusp,
$Y^{loc}$ and $Z^{loc}$
are both of type $A$
and  $Y^{loc}\cap Z^{loc}$ is nodal.

For \DHrefpart{iii} it is enough to show that $\sX''_0$ is Gorenstein.
This follows from \DHrefpart{i}.
\end{proof}
\end{lemma}

\begin{lemma}\label{9.6}\label{9.6-}
\part[i] $\sX'$ is Gorenstein,
$Y'$ and $Z'$ have only singularities of type $A$
and their canonical classes are linearly equivalent to zero.
\part[ii] 
The fibre $f'^{-1}(P')$
is either an elliptic curve or a cycle of rational curves.
\begin{proof}
\DHrefpart{i}: By Lemma \ref{previous}
$\sX''_0$ has Gorenstein singularities so $\sX'_0$ is
everywhere locally Gorenstein. Since $Y'$ and $Z'$
have singularities of type $A$ so do $Y'$ and $Z'$.

Now $f':\sX'\to S'$ is generically Jacobian,
and so Jacobian over $S'-\{P'\}$.
Therefore $K_{\sX'-f'^{-1}(P')}$ is linearly equivalent
to the pull back of a line bundle on $S'-\{P'\}$,
and so $K_{\sX'-f'^{-1}(P')}\sim 0$. Since $\sX'$
is Gorenstein, $K_{\sX'}\sim 0$; the triviality of
$K_{Y'}$ and $K_{Z'}$ is then immediate.

\DHrefpart{ii}: The fact that $f'^{-1}(P')$ is reduced and nodal 
follows from the description given in
Lemma \ref{previous}. The rest follows from the triviality of $K_{\sX'}$.
\end{proof}
\end{lemma}

\begin{lemma}\label{9.8+1}
  \part[i] $f':\sX'\to S'$ is a semi-stable family
  of curves of genus $1$.
  \part[ii] $f^h:\sX^h\to S^h$ is the geometric quotient of $f':\sX'\to
  S'$
  by an equivariant action of a cyclic group of order $1,2,3,4$ or $6$.
  
  \part[iii] Suppose that $Y\to C$ is
  generically smooth and has finite local
  monodromy around $P$. Then $f'^{-1}(P')$ is smooth
  and $Z\to D$ is
  generically smooth with  finite local
  monodromy around $P$.
  \part[iv] If $f'^{-1}(P')$ is smooth, then $Y\to C$ and $Z\to D$
  are both generically smooth with finite local monodromy around $P$
  and the fibres of the pair $(Y\to C,\ Z\to D)$ over $P$
  are of types either $(I_0,I_0)$
  or $(\bD_4,\bD_4)$ or $(\bR,\bR^*)$, where $R=II,III,IV$.
  \begin{proof} This is an immediate consequence of Lemma
    \ref{9.6} (2).
  \end{proof}
\end{lemma}

\begin{lemma}\label{9.8+2}
  If $Y$ is birational to a special elliptic surface,
  then $f'^{-1}(P')$ is smooth and the fibres of the pair
  $(Y,Z)$ over $P$ are of type either $(\bD_4,\bD_4)$
  or $(I_0,I_0)$.
  \begin{proof} This follows at once from Lemma \ref{9.8+1}.
  \end{proof}
\end{lemma}

Now suppose that, under the period map, the image of $0\in\Delta$
under the period map is the direct sum of matrices of special elliptic
surfaces $\tV_1,...,\tV_r$.

\begin{lemma} If $f_i:X_i\to C_i$
  is generically smooth and $h^2(X_i,\sO_{X_i})>0$
  then $X_i$ is a special elliptic surface.
  \begin{proof} The Hodge structure on $H^2(X_i)$
    embeds into the direct sum $\oplus H^2(\tV_i)$.
    Since each $\tV_i$ is special, it follows that the
    global monodromy on $H^1$ of the generic fibre
    of $f_i:X_i\to C_i$ is a group of order $2$,
    so that $X_i$ is a special elliptic surface. (Recall that
    $f_i:X_i\to C_i$ has a section, so that if $X_i$ is not simply
    connected
    then it is ruled, which would contradict $h^2(X_i,\sO_{X_i})>0$.
    So $X_i$ is simply connected.)
  \end{proof}
\end{lemma}

\begin{lemma} If $X_i$ is birational to a special elliptic surface
  $U$ then it is isomorphic to the $\bD_4$-model of $U$.
  \begin{proof}
    This follows from Lemma \ref{9.11} and Lemma \ref{9.8+2}.
  \end{proof}
\end{lemma}

We now prove Theorem \ref{10.x}.

Say $\sS=\{i:\ f_i\ \textrm{is\ generically\ smooth\ and}\ h^2(X_i,\sO_{X_i})>0\}$.
So $\sS$ is the set of indices $i$ such that $X_i$ is not a union of
ruled surfaces. Say $\vert\sS\vert =s$.

Recall that $K_{\sX/\Delta}\vert_{X_i}\sim f_i^*(\a_i)$,
where $\a_i\in\Q_{>0}$. Suppose $i\in\sS$; then the fibres
of $f_i:X_i\to C_i$ are all of type either $\bD_4$ or $I_0$.
Suppose that $X_i$ has $b_i$ double curves of type $\bD_4$
and $c_i$ of type $I_0$.
Then, by the adjunction formula,
$$K_{X_i}\sim f_i^*(\a_i-{\half}b_i-c_i).$$
So $h_i-1=\a_i-{\half}b_i-c_i.$

The objects $f_j^{-1}(C_j)$ form a connected tree,
so, for every pair $i,j\in\sS$, there is a unique path
that links them. So there is a total of
$2\times (s-1)$ double curves on $\sqcup_{i\in\sS}X_i$
each of which is linked to a
double curve on some other surface $X_j$
with $j\in\sS$. Therefore
\begin{eqnarray*}
  \sum_k\a_k&=&h-1=\sum_{i\in\sS} h_i-1=\sum_{i\in\sS}(h_i-1)+(s-1)\\
  &=& \sum_{i\in\sS}(\a_i-{\half}b_i-c_i)+(s-1).
\end{eqnarray*}
It follows that
$${\half}\sum_{i\in\sS}b_i+\sum_{i\in\sS}c_i+\sum_{k\not\in\sS}\a_k=s-1.$$
Since $\sum b_i+\sum c_i$ is the total number of
double curves on $\sqcup_{i\in\sS}X_i$, we have
$$\sum b_i+\sum c_i\ge 2(s-1).$$
Comparison of these shows that $c_i=0$ for all $i\in\sS$
and, since every $\a_k>0$, that
$\sS$ is the full set of indices. That is,
$\sX_0$ is the union of the special elliptic surfaces $X_1,...,X_s$. 

  Therefore the direct sums $\oplus_1^s H^2_{prim}(X_i)$
  and $\oplus_1^rH^2_{prim}(\tV_i)$ of integral Hodge structures
  are isomorphic. Moreover, there are elliptic curves $E,F$
  and hyperelliptic curves $A_1,...,A_s,B_1,...,B_r$ such that
  $X_i\cong[(E\times A_i)/\langle\iota\rangle]$
  and $V_j\cong[(F\times B_j)/\langle\iota\rangle]$.
  The product $E\times A_i$ can then be recovered from $X_i$
  as the unique double cover of $X_i$ ramified in its nodes,
  and the same holds for $F\times B_j$. Therefore
  there is an isomorphism
  $$H^1(E)\otimes\left(\oplus H^1(A_i)\right)\to
  H^1(F)\otimes\left(\oplus H^1(B_j)\right)$$
  of integral Hodge structures of weight $2$.

   From Theorem \ref{taylor1},
   to be proved in Section \ref{taylor} below, and the Torelli
   theorem for curves, it follows that
  $E\cong F$, that
  $r=s$ and that, after re-ordering if necessary,
  $A_i\cong B_i$ for every $i$. So $X_i$ is birational to $V_i$,
  and Theorem \ref{10.x} is proved.
\end{proof}
\end{theorem}

We also prove some further results about the
object $\sY$ of Theorem \ref{10.x} that will be useful later.

Say that $V_i$ meets $s_i$ other surfaces $V_j$ in $\sY_0$.
That is, $V_i$ contains $s_i$ double curves $\delta_{ij}$, 
each of which is of type $\bD_4$. 
Since each $V_i$ has only nodes, we can assume that,
after some finite base change $\Delta\to\Delta$ if necessary,
each $V_i$ is smooth outside
the double curves $\delta_{ij}$.

\begin{lemma} If $B_i\subset V_i$ is a section
of $f_i:V_i\to C_i$ then $B_i$ meets each
$\delta_{ij}$ in a node and $B_i^2=-(h_i+1)+s_i/2$.
\begin{proof} The strict transform $\tB_i$ of $B_i$ on $\tV_i$
is a section, so meets each $\tD_4$ fibre in an end curve
of the $\tD_4$ configuration. So $B_i$ meets
$\delta_{ij}$ is a node. Since $\tB_i^2=-(h_i+1)$
the lemma follows.
\end{proof}
\end{lemma}

\begin{proposition}\label{10.y} $\sY\to T$ has a section.
  \begin{proof} There is a section of $\sX_\eta$, and so
    an irreducible Weil divisor $D\subset\sY$
that restricts to a section of $\sY\to T$
over the generic point of $\Delta$.
We can write $D\cap V_i=D_i+\psi_i$
where $D_i$ is a section of 
$f_i:V_i\to C_i$ and $\psi_i$ is
an $f_i$-vertical curve that meets $D_i$.
To show that $D$ is a section it is enough to show that
each $\psi_i=0$.

Recall that $p_g(V_i)=h_i$ and $p_g(\sY_t)=h$ for $t\ne 0$.
Then $D^2.\sY_t=-(h+1)$
and $(D_i)^2_{V_i}=-(h_i+1)+s_i/2$,
since $D_i$ passes through $s_i$ nodes on $V_i$.
So
\begin{eqnarray*}
-(h+1)&=&D^2.\sY_t=\sum D^2.V_i=\sum (D_i+\psi_i)^2_{V_i}\\
&=&-\sum (h_i+1)+\sum s_i/2+2\sum(D_i.\psi_i)_{V_i}.
\end{eqnarray*}
Now $\sum s_i=2(r-1)$ and $\sum h_i=h$, so that
$\sum(D_i.\psi_i)_{V_i}=0$. It follows that each $\psi_i=0$.
\end{proof}
\end{proposition}
\begin{lemma}\label{10.z}
The singularities of $\sY$ are of index at most $2$.
\begin{proof} Inspection.
\end{proof}
\end{lemma} 
\end{section}
\begin{section}{Unscrewing tensor products of weight $1$ Hodge structures}\label{taylor}
  Here we prove a result that yields a Torelli theorem for sesqui-special elliptic surfaces.
The phrases  ``Principally polarized abelian variety'' and
  ``irreducible representation''
  are abbreviated
  to ``ppav''
  and ``irrep''.

  \begin{theorem}\label{taylor1}
    Suppose that $A_0,B_0$ are elliptic curves, that
    $A_1,...,A_r,B_1,...,B_s$
    are ppav's and that $\Hom(A_i,A_j)=\Z\delta_{ij}=\Hom(B_i,B_j)$
    for all $i,j\ge 0$. Say $A=\sum_{i>0}A_i$ and $B=\sum_{j>0}B_j$
    and give them the product principal polarizations.
    Assume that $\dim A=h\ge 2$ and that there is an isomorphism
    $H^1(A_0)\otimes H^1(A)\to H^1(B_0)\otimes H^1(B)$
    of ${\Z}$-Hodge structures.

    Then $A_0$ is isomorphic to $B_0$ and $A$ is isomorphic to $B$
    as ppav's.

    \begin{proof} We first prove a result (Proposition \ref{Proposition 1}) about polarizable
      $\Q$-Hodge structures, and then at the end clear denominators.

      We shall use some elementary facts about Mumford--Tate groups
      of polarizable Hodge structures, as described in \cite{Mo}.

      The category ${\Q}{\rm{HS}}^{pol}$ of polarizable ${\Q}$-Hodge structures
is Tannakian and semi-simple; the corresponding group $MT$
(the universal Mumford--Tate group) is a pro-reductive group over ${\Q}$.
There is a cocharacter $\nu: {\GG}_{m,\Q}\to MT$ such that if $\rho_V:MT\to GL(V)$
is the representation corresponding to the polarizable ${\Q}$-Hodge structure $V$
of pure weight $n$
then $\rho_V(\nu(r))(v)=r^{-n}v$ for all $v\in V$. The image of $\rho_V$ is $MT(V)$.

Moreover the Deligne torus ${\SS}=R_{\C/\R}\GG_m$
embeds into $MT_{\R}$ so as to extend
the embedding $\nu_{\R}$ of $ {\GG}_{m,\R}$. The action
of $ {\SS}_{\C}$ on the subspace $V^{p,q}$ of $V_{\C}$
is given by $(z_1,z_2)(v)=z_1^{-p}z_2^{-q}v$ and ${\GG}_{m,\C}$
appears as the diagonal
subgroup of $ {\SS}_ {\C}$.

We have polarized ${\Z}$-Hodge structures $A_i'=H^1(A_i,\Z)$, $B_j'=H^1(B_j,\Z)$
of weight $1$ and dimensions $2,2,2h,2h$.
Tensor these Hodge structures with ${\Q}$; write
$A_i=A_i'\otimes{\Q}$, etc., so that there are decompositions
$$A=\sum_1^r A_i\ {\rm{and}}\ B=\sum_1^s B_j$$
into irreducible components.
Assume that $A_0\otimes A$ is isomorphic to $B_0\otimes B$
as ${\Q}$-Hodge structures,
that $h\ge 2$ and that $MT(A_0)$ and $MT(B_0)$ are both isomorphic to
$GL_2$.

\begin{proposition}\label{Proposition 1}
$A_0$ is isomorphic to $B_0$ and $A$ is isomorphic to $B$ as
${\Q}$-Hodge structures.
\begin{proof}
  There are representations
$\rho=\rho_{A_0}:MT\to GL(A_0)$ and $\sigma_i=\rho_{A_i}:MT\to GL(A_i)$.
By assumption, $\rho$ is surjective.
Say $G_i=\im(\sigma_i)=MT(A_i)$ and
$H_i=\im(\rho\oplus\sigma_i) =MT/(\ker(\rho)\cap \ker(\sigma_i)=MT(A_0+A_i)$.

Consider $\rho\otimes\sigma_i:H_i\to GL(A_0\otimes A_i)$ and let $Z_i\subset H_i$ denote its kernel.

\begin{lemma} $Z_i$ is central in $H_i$.
  \begin{proof}
    It is enough to show that, if
    $x$ and $y$ are two matrices such that $x \otimes y = 1$, then $x$
    and $y$ are scalars.
    
Every eigenvalue of $y$ must be the inverse of every eigenvalue of $x$,
so $x$ and $y$ each have a unique eigenvalue, which we can assume is $1$. 
If $x \ne 1$ consider a vector $u$ with
$(x-1)u\ne 0$ and a non-zero vector $v$ with $yv=v$. Then
$(x \otimes y) (u \otimes v)\ne u \otimes v$, a contradiction.
\end{proof}
\end{lemma}

\begin{proposition}\label{well known}
  \part[i] The derived subgroup $G^{der}$ of a reductive group $G$
  defined over a field has a simply connected universal cover $G^{sc}$.

   \part[ii]  Any simply connected (or adjoint) reductive group over a
   field is a product of
simple simply connected (or adjoint) reductive groups. 

\part[iii]  If $G_1,G_2$ are groups  over a field
then every irrep $W$ of $G_1\times G_2$ is of the form
$W=V_1\otimes V_2$
where $V_i$ is an irrep of $G_i$.
\begin{proof} \DHrefpart{i}: For split groups this is pointed out in 7.3.4 of \cite{SGA3}
  XXIV. In general see Corollary A.4.11 of \cite{CGP}.
  
\DHrefpart{ii}: For adjoint groups this is
Prop. 5.5 of \cite{SGA3} XXIV. As stated in 5.3 of \emph{op. cit.}
the same holds for simply connected groups.

\DHrefpart{iii}: This is well known.
\end{proof}
\end{proposition}
  
The construction of $G^{sc}$ is functorial in $G$, so $MT^{sc}$ exists, and the homomorphisms
$$MT\to H_i\to GL(A_0)\times MT(A_i)$$
give homomorphisms
$$MT^{sc}\to H_i^{sc}\buildrel{\pi}\over\to MT(A_0)^{sc}\times MT(A_i)^{sc}=SL_2\times MT(A_i)^{sc}$$
of simply connected reductive groups over ${\Q}$,
and $H_i^{sc}$ maps surjectively to both factors.

\begin{lemma}\label{Claim 3}
  If there is an $SL_2$-factor in $MT(A_i)^{sc}$ then $MT(A_i)^{sc}=SL_2$
and $\dim(A_i)=2$.
\begin{proof} Suppose $MT(A_i)^{sc}=SL_2\times K_i$, where $K_i$ is simply connected.
  Since $\End(A_i)=\Q$, it follows that the connected part $Z_i$ of the
  centre of $MT(A_i)$ is $Z_i\cong\GG_m$. Therefore $MT(A_i)/MT(A_i)^{der}\cong\GG_{m}$,
  so that $A_i$ is also an irrep of $MT(A_i)^{sc}$.

  There are isogenies
  $$Z_i\times SL_2\times K_i\to Z_i\times MT(A_i)^{der}\to MT(A_i),$$
  so that
  $A_i=X_i\otimes Y_i$
where $X_i$ is an irrep of $Z_i\times SL_2$ and $Y_i$ is an irrep of $K_i$.

Consider also the isogeny
$${\widetilde{\SS}}=\SS\times_{MT(A_i)_\R}(Z_i\times SL_2)_\R
\times K_{i,\R} \to \SS.$$
Say that ${\widetilde{\SS}}_{\C}$ acts on $X_{i,\C}$
with distinct characters $\chi_1,...,\chi_a$
and on $Y_{i,\C}$
with distinct characters $\psi_1,...,\psi_b$.
Since ${{\SS}}_{\C}$ acts on $A_{i,\C}$
via just 2 characters,
$\alpha_1:(z_1,z_2)\mapsto z_1^{-1}$
and $\alpha_2:(z_1,z_2)\mapsto z_2^{-1}$, it follows that
$\{\chi_i+\psi_j\}=\{\alpha_1,\alpha_2\}$,
so that $a,b\le 2$.


It is easy to see that $a=b=2$ is impossible, so that
there are two possibilities.

\begin{enumerate}
\item $a=1$. Then there is a non-trivial torus in
  $(\GG_m\times SL_2)_\C$
  that acts on $X_{i,\C}$
  with only $1$ weight. But $X_i$
  is an irrep of $\GG_m\times SL_2$, so it is trivial.
But this is a contradiction
since $X_i\otimes Y_i$ is a faithful rep of $MT(A_i)$.

\item $b=1$. The same argument shows that $Y_i\otimes{\C}=Y_i^{p,p}$
for some $p$, so that $Y_i(p)$ is trivial, so $1$-dimensional. Then $K_i=1$,
so $MT(A_i)^{sc}=SL_2$ and $X_i=A_i$ is an irrep of $SL_2$ for which some
non-trivial torus in $SL_2$ has only two weights. So $\dim(A_i)=2$
and Lemma \ref{Claim 3} is proved.
\end{enumerate}
\end{proof}
\end{lemma}

\begin{lemma} The homomorphism $\pi:H_i^{sc}\to MT(A_0)^{der}\times MT(A_i)^{sc}$
is an isomorphism and $A_0\otimes A_i$ is irreducible.
\begin{proof}  Since both projections $H_i^{sc}\to MT(A_0)^{der}$
and $H_i^{sc}\to MT(A_i)^{sc}$ is surjective, it follows from
the claim that we have just established
that {\it{either}} $\pi:H_i^{sc}\to MT(A_0)^{der}\times MT(A_i)^{sc}$ is an isomorphism
{\it{or}} $G_i^{sc}=SL_2$ and the image of $\pi$ is the diagonal
subgroup of $SL_2\times SL_2$. In this latter case, however, we have
$A_0\cong A_i$, contradiction.
\end{proof}
\end{lemma}

So the equations
$$A_0\otimes A=\sum_1^r A_0\otimes A_i\ \textrm{and}\ B_0\otimes B=\sum_1^s B_0\otimes B_j$$
are decompositions into irreducibles.
So $r=s$ and $A_0\otimes A_i\cong B_0\otimes B_i$
(after re-ordering if necessary).

Since $A_0\otimes A_i$ is irreducible,
we can recover the factors $A_0$ and $A_i$ as representations
of $MT$ up to tensoring by mutually inverse characters. 
Considering the action of ${\SS}$ determines these characters.
Since also $\dim(A)=\dim(B)>2$,
it follows that 
that $A_i\cong B_i$ for all $i\ge 0$.

This completes  the proof of Proposition \ref{Proposition 1}.
\end{proof}
\end{proposition}

\begin{lemma} The terms $A_0\otimes A_i$ are pairwise non-isomorphic.
  \begin{proof}
The image of $MT(A_0+A_1+A_2)^{der}$ in $MT(A_0)\times MT(A_1)\times MT(A_2)$
equals $MT(A_0)^{der}\times MT(A_1)^{der}\times MT(A_2)^{der}$,
since $\Hom(A_i,A_j)=0$ for $i\ne j$.
The result follows.
\end{proof}
\end{lemma}

\begin{lemma}
  \part[i] $\End(A_0\otimes A_i)={\Q}$ for $i>0$ and
  $\End(\sum_{i>0} A_0\otimes A_i) =\End(A)={\Q}^r$.

  \part[ii] Every homomorphism $\phi:A_0\otimes A\to B_0\otimes B$
of ${\Q}$-Hodge structures 
can be written as $\phi=\phi_1\otimes\phi_2$
where $\phi_1:A_0\to B_0$ and $\phi_2:A\to B$
are homomorphisms of ${\Q}$-Hodge structures.
\begin{proof} \DHrefpart{i} follows at one from
the equality $MT(A_0\otimes A_i)^{der}=MT(A_0)^{der}\times MT(A_i)^{der}$
and the fact that
$\Hom(A_0\otimes A_i,A_0\otimes A_j)=\Hom(A_i,A_j)={\Q}\delta_{ij}$.
\DHrefpart{ii} is an immediate consequence.
\end{proof}
\end{lemma}

\begin{proposition}
Suppose that $\phi':A_0'\otimes A'\to B_0'\otimes B'$
is an isomorphism of ${\Z}$-Hodge structures.
Then there are isomorphisms $\phi_1':A_0'\to B_0'$
and $\phi_2':A'\to B'$ of ${\Z}$-Hodge structures
such that $\phi'=\phi_1'\otimes\phi_2'$.
\begin{proof}
  Write $\phi=\phi'\otimes{\Q}$
and $\phi=\phi_1\otimes\phi_2$.
The structure theorem for modules over
a PID shows that there are ${\Z}$-bases
$\{e_i\}$, $\{f_i\}$, $\{a_j\}$ and $\{b_j\}$
of $A_0',B_0',A'$ and $B'$ and scalars $\lambda_i,\mu_j\in {\Q}^*$
such that $\phi_1(e_i)=\lambda_if_i$
and $\phi_2(a_j)=\mu_jb_j$
for all $i,j$. 
Then $\phi'(e_i\otimes a_j)=\lambda_i\mu_j f_i\otimes b_j$,
so that $\lambda_i\mu_j=\pm 1$ for all $i,j$. It follows that
$\lambda_i=\pm \lambda_1$ and $\mu_j=\pm \mu_1$
for all $i,j$. After replacing $f_i$ by $\pm f_i$
and $b_j$ by $\pm b_j$ if necessary, we see
that there are scalars $\lambda,\mu$ such that
$\phi_1(e_i)=\lambda f_i$,
$\phi_2(a_j)=\mu b_j$ for all $i,j$
and $\lambda\mu=1$.

Then $\phi'=\phi_1'\otimes\phi_2'$ where
$\phi_1'=\lambda^{-1}\phi_1$ and
$\phi_2'=\lambda\phi_2$,
and each $\phi'_i$ is a ${\Z}$-isomorphism.
\end{proof}
\end{proposition}

Now we can complete the proof of the theorem. We have ppav's
$A_0,A_1,...,A_r$ and $B_0,B_1,...,B_s$ such that
$\Hom(A_i,A_j)=\Z\delta_{ij}=\Hom(B_i,B_j)$.
There is an isomorphism
$H^1(A_0)\otimes\sum H^1(A_i)\to H^1(B_0)\otimes\sum H^1(B_j)$
of ${\Z}$-Hodge structures and $\dim\sum_{i>0} A_i\ge 2$.
We want to prove that
$A_0\cong B_0$ and $\sum_{i>0} A_i$ is isomorphic to $\sum_{j>0} B_j$
as ppav's.

For this, it only remains to observe that, if $A=\sum A_i$
where each $A_i$ is a ppav
and $\Hom(A_i,A_j)=\Z\delta_{ij}$ for all $i,j$, then
the only principal polarization on $A$ is the standard
product principal polarization.
\end{proof}
\end{theorem}
\end{section}
\begin{section}{The main results for surfaces}\label{consequences}
Our goal is to give a first order description of the period locus $PL_h$ 
in a neighbourhood of $\sW_{1^h}$, or
of $\sW_{1^h}/(SL_2(\Z)\wr\Symm_h)\times SL_2(\Z)$.

\begin{definition}\label{RDP, no monodromy}
Denote by
$\sJE^{RDP,n.m.}$ (resp., $\sJE^{sm,n.m.}$) the stack
of \emph{generalized} RDP (resp., smooth) semi-stable
Jacobian elliptic surfaces with no monodromy.
Its objects over a scheme $\Delta$
are pairs $(\sX{\buildrel{f}\over{\to}} S\to\Delta,S')$ 
with the following properties:
\begin{enumerate}
\item $\sX\to S$ and $S\to\Delta$ are projective;

\item $\sX\to\Delta$ and $S\to\Delta$ are flat, Cohen--Macaulay and
of relative dimensions $2$ and $1$, respectively;

\item $S'$ is a section of $\sX\to S$;

\item $S\to\Delta$ is a semi-stable family of curves of compact type;

\item for every geometric point $\delta$ of $\Delta$
the fibre $\sX_\delta$ of $\sX\to\Delta$
is a reduced sum $\sX_\delta=\sum V_i$ of irreducible components $V_i$
each of which
is of the form $V_i=f^{-1}(C_i)$ for some irreducible component
$C_i$ of $S_\delta$;

\item over the complement of the
nodes of $S_\delta$ each
$V_i$ has only RDPs (resp., is smooth)
and $V_i\to C_i$ is relatively minimal;

\item $2S'$ is a Cartier divisor on $\sX$;

\item $2S'$ is $f$-ample (resp., is $f$-nef);

\item the restriction $f_i=f\vert_{V_i}:V_i\to C_i$ 
gives the minimal resolution $\tV_i$ of $V_i$ 
the structure of an elliptic surface of genus $h_i$ over $C_i$
on which the strict transform of $S'\cap V_i$ is the identity section;

\item $\sX_\delta$ has slc singularities;

\item if $P=C_i\cap C_j$ is a node of $S_\delta$ then $f^{-1}(P)$
is a fibre of type $\bD_4$ on each of $V_i$ and $V_j$
and the section $S'$ meets $f^{-1}(P)$ in one of the four points on it
that is an $A_1$ singularity on both $V_i$ and $V_j$;

\item the rank one sheaf $\omega_{\sX/\Delta}^{[2]}=\sO(2K_{\sX/\Delta})$ is
  the pullback of an
  ample
  invertible sheaf on $S$;

\item $\sum h_i$ is independent of the point $\delta\in\Delta$.
  (This last is the ``no monodromy'' condition.)
\end{enumerate}
\end{definition}

Observe that there is a commutative diagram
$$\xymatrix{
  {\sJE^{sm}}\ar@<-1ex>[d]\ar@{^{(}->}[r]^-{open}&
  {\sJE^{sm,n.m.}}\ar@<-3ex>[d]\\
  {\sJE^{RDP}}\ar@{^{(}->}[r]^-{open}&
  {\sJE^{RDP,n.m.}}
}$$
where each vertical arrow is given by passing to the relative canonical model.
It that describes the fact that $\sJE^{RDP,n.m.}$ (resp., $\sJE^{sm, n.m.}$)
is a partial compactification of the stack
$\sJE^{RDP}$ (resp., $\sJE^{sm}$).

\begin{proposition}\label{extend period}
  \part[i]
  The relative canonical model $\sY\to T\to\Delta$ of
  the object $\sY\to\Delta$ of Theorem \ref{10.x}
  (an output of the MMP)
  is an object over $\Delta$ of $\sJE^{RDP,n.m.}$.

  \part[ii] 
  The relative canonical model $\bsY\to \sB\to\Delta$ of
  the object $\bsY\to\Delta$ of Corollary \ref{8.12}
  (an output of a plumbing construction)
  is an object over $\Delta$ of $\sJE^{RDP,n.m.}$
  (resp., $\sJE^{sm,n.m.}$)
  if, away from the curve $\bY_a\cap\bY_b=\overline{\phi}$,
  each surface $\bY_i$ is an RDP (resp., smooth) elliptic surface.
  \begin{proof} For $\bsY$ of Corollary \ref{8.12}
    this is clear. For $\sY$ this follows from Theorem \ref{10.x},
    Proposition \ref{10.y} and Lemma \ref{10.z}.
  \end{proof}
\end{proposition}
\begin{proposition}\label{10.1} $\sJE^{RDP,n.m.}$ is smooth.
\begin{proof} It is enough to find a smooth stack $\sS$ and
a surjective smooth morphism $\sS\to\sJE^{RDP,n.m.}.$

Recall that a stable hyperelliptic curve of compact type 
over a base $\Delta$
is a stable curve $q:C\to\Delta$ of compact type with a $q$-equivariant
involution $\iota$ 
such that $\iota$ covers the identity map $id_\Delta$,
the fixed point locus of $\iota$ is finite over $\Delta$,
the geometric quotient $[C/\iota]\to\Delta$ is a pre-stable curve of genus zero
and, for every geometric point $\delta\in\Delta$, 
$\iota$ preserves each component of the geometric fibre $C_\delta$.
In particular, every node of $C_\delta$ is, therefore, $\iota$-invariant.
These are the objects of the stack $\sHyp^c=\coprod_g\sHyp_g^c$, which
is smooth and
is the closure of the hyperelliptic locus in
the stack of stable curves of compact type.

Consider the stack $\sS$ whose objects over $\Delta$ are triples 
$(Y{\buildrel{p}\over{\to}}C{\buildrel{q}\over{\to}}\Delta,\iota,\iota_Y)$ 
with the properties that
\begin{enumerate}
\item
$(q:C\to\Delta,\iota)$ is a stable hyperelliptic curve of compact type,

\item $p:Y\to C$ is flat, projective and of relative dimension $1$,

\item there is a section $C_0$ of $p$ 
that lies in the relatively smooth locus of $p$
and is a $p$-ample Cartier divisor on $Y$,
 
\item $p$ is smooth over the fixed points of $\iota$
(this includes all nodes on each $C_\delta$),

\item for each geometric point $\delta$ of $\Delta$
and each irreducible component $D$ of $C_\delta$,
the inverse image $p^{-1}(D)$ has only RDPs,

\item for each $D$ as just described, $p^{-1}(D)\to D$
is a Jacobian elliptic surface on which $C_0$
is the identity,

\item $\iota_Y$ is an involution of $Y$ that lifts $\iota$
and preserves $C_0$ and

\item for every geometric point $\delta\in\Delta$
and for each node $x$ in $C_\delta$,
$\iota_Y$ acts on the fibre $p^{-1}(x)$ as $(-1)$.
\end{enumerate}
It follows that $\iota_Y$ acts as $(-1)$ on every smooth fibre of $p$
and that the locus of fixed points of $\iota_Y$ consists
of the $2$-torsion points on the fibres $p^{-1}(x)$ over the fixed
points $x$ of $\iota$.
\begin{lemma}
The forgetful morphism $\sS\to\sHyp^c$ is smooth.
\begin{proof}
$\sS\to \sHyp^c$ is the base change under $\sHyp^c\to\sM^c$
of the morphism $\sJE^{RDP,sst}\to\sM^c$ from
Lemma \ref{sm via sextic}. The result now follows from that lemma.
\end{proof}
\end{lemma}
Since $\sHyp^c$ is smooth the same is true
of $\sS$. 

Taking geometric quotients by $\iota$ gives a morphism
$\pi:\sS\to\sJE^{RDP,n.m.}.$

\begin{lemma}
$\pi:\sS\to\sJE^{RDP,n.m.}$ is smooth.
\begin{proof}
Given an object $\sX{\buildrel{f}\over{\to}} S\to\Delta$ in
$\sJE^{RDP,n.m.}$
the objects in its pre-image
are obtained by taking a double cover of $\sX$ whose branch locus
is an $f$-vertical divisor $B$ and a finite scheme that is disjoint
from $B$ and supported on the $A_1$ singularities in the $\bD_4$-fibres
of $f$. 

Suppose that $\Delta\inj\Delta_1$ is a thickening
of Artin schemes and that we are given
$\sX_1{\buildrel{f_1}\over{\to}} S_1\to\Delta_1$
that restricts to $\sX{\buildrel{f}\over{\to}} S\to\Delta$.
Suppose also that we have an object $Y$ in the pre-image
of $\sX{\buildrel{f}\over{\to}} S\to\Delta$ as above.
Certainly
$B$ can be lifted to a vertical divisor
$B_1$, so that $Y$ can be extended
to a double cover $Y_1^0$ over $\Delta_1$ outside
the singular points on the $\bD_4$-fibres.
Moreover, the local triviality of $\omega_{\sX_1/S_1}^{[2]}$
implies that $Y_1^0$ can be extended across these singular points.
\end{proof}
\end{lemma}
This completes the proof of Proposition \ref{10.1}.
\end{proof}
\end{proposition}
Now restrict to the closed substacks $\sJE^{*,n.m.}_{s.c.}$
of $\sJE^{*,n.m.}$ consisting
of simply connected Jacobian elliptic surfaces,
where $*$ denotes $sm$ or $RDP$,
and write
$$\sJE^{*,n.m.}_{s.c.}=\coprod_h\sJE^{*,n.m.}_h$$
where $h$ denotes the geometric genus.
Then the $2$-Cartesian diagram above restricts to give
a $2$-Cartesian diagram
$$\xymatrix{
  {\sJE^{sm}_h}\ar@<-1ex>[d]\ar@{^{(}->}[r]^-{open}&{\sJE^{sm,n.m.}_h}\ar@<-3ex>[d]\\
{\sJE^{RDP}_h}\ar@{^{(}->}[r]^-{open}&{\sJE^{RDP,n.m.}_h}
}$$

Recall the loci $\sW_{h_1,...,h_r}$ described in Definition \ref{defn of sW}.
Recall also that $h=\sum h_i$.

\begin{proposition}\label{11.a}
The period map $per:\sJE^{sm}_h\to\sV_h/\mathfrak G$
extends to a morphism
$per^+:\sJE^{sm,n.m.}_h\to\sV_h/\mathfrak G$
that is proper over a neighbourhood of
the generic point of
each $\sW_{h_1,...,h_r}$
and that fits into a $2$-commutative diagram
$$\xymatrix{
{\sJE^{sm,n.m.}_h}\ar[r]^-{per^+}\ar@<-3ex>[d]&{\sV_h/\mathfrak G}\ar[d]\\
{\sJE^{RDP,n.m.}_h}\ar[r]&{[\sV_h/\mathfrak G].}
}$$ 
\begin{proof} The properness follows from Theorem \ref{10.x}
  and Proposition \ref{extend period}.
\end{proof}
\end{proposition} 
\begin{lemma}\label{sm is smooth}
For each $(h_1,...,h_r)$
and each point $x$ in $\sW_{h_1,...,h_r}$,
the stack
$\sJE^{sm,n.m.}_h$ is smooth at each point lying over
$x$.
\begin{proof}
Suppose that $(X\to S,S')$ is a geometric point of 
$\sJE^{sm,n.m.}_h$ 
that maps to $x$.
Let $r:X\to Y$ be the contraction of all vertical
$(-2)$-curves that lie in the smooth locus
of $X$ and are disjoint from $S'$.
Recall the combinatorial description of $Y$.

\begin{enumerate}
\item
$X=\sum_1^r X_i$ and
$Y=\sum_1^r Y_i$ where, by assumption, each $Y_i$ is
birational to a geometric quotient
$[(E_i\times C_i)/\iota]$. 

\item The elliptic curves $E_i$
are all isomorphic since $Y$ is connected.

\item The configurations $X=\sum X_i$, $Y=\sum Y_i$
and $S=\sum C_i$
are isomorphic trees.

\item $S''=r(S')$ is a section
of $Y\to S$. $2S''$ is Cartier
and is ample relative to $S$.

\item Say $\s_{ij}=X_i\cap X_j$ if this is non-empty
and $\phi_{ij}=r(\s_{ij})$.
Then $\s_{ij}$ and $\phi_{ij}$ are fibres of type $\bD_4$
on each of $X_i,X_j, Y_i,Y_j$, as appropriate. 

\item $Y_i$ has $4$ singularities of type $A_1$ on each
$\phi_{ij}$ and has $D_4$-singularities disjoint
from $S''$ and from the $\phi_{ij}$.
\end{enumerate}

The divisor $2S''$ defines
a finite morphism $\rho:Y\to Z=\sum_1^r Z_i$ of degree $2$.
Say $\psi_{ij}=\rho(\phi_{ij})$. Then, 
because of the nature of special Jacobian elliptic surfaces,
\begin{enumerate}
\item $Z\to S$ is a $\P^1$-bundle,

\item the branch locus $B\subset Z$
is $B=B_0+B_1+\sum_{ij}\psi_{ij}$ where

\item $B_0=\rho(S'')$ is
a section of $Z\to S$
and a Cartier divisor on $Z$, 

\item $B_1$ is disjoint from $B_0$,

\item $B_1$ is a sum $B_1=\sum_1^3 D_i$
of three sections $D_i$ of $Z\to S$,

\item each $D_i$ is a Cartier divisor on $Z$,

\item the $D_i$ are linearly equivalent and

\item $B_1$ has ordinary triple points over the complement of 
the nodes in $S$ and, as does $B_0$, 
meets each  $\psi_{ij}$ transversely.
\end{enumerate}

Conversely, given $(Z,B)$, we recover $Y$ as the double cover
of $Z$ branched along $B+\sum_{ij}\psi_{ij}$.

We next prove that we can deform the triple points of
$B_1$ independently, while fixing $Z$. 
For this, let $\Sigma$ denote the set of
triple points of $B_1$ and $\sI_\Sigma$ its sheaf of ideals, 
and consider the short exact sequence
$$0\to\sI_\Sigma^3(B_1)\to\sO_Z(B_1)\to\sO_Z/\sI_\Sigma^3\to 0$$
of coherent sheaves on $Z$.
It is straightforward to verify that
$H^0(Z,\sI_\Sigma^3(B_1))=\Symm^3 H^0(Z,\sI_\Sigma(D_1))$
and that then a count of dimensions shows that the map
$H^0(Z,\sO_Z(B_1))\to H^0(Z,\sO_Z/\sI_\Sigma^3)$ is surjective.

It follows that the morphism 
$\Def_Y\to\prod\Def_{Y,P}$
of deformation spaces,
where $P$ runs over the $D_4$ singularities of $Y$,
is formally smooth.
\end{proof}
\end{lemma}

Now fix an alkane $\G$ of genus $h$.
Recall the period space $\sV_h$ and its $9h+9$-dimensional closed subvariety
$\sV_\Gamma$ whose points correspond to configurations of
special elliptic Kummer surfaces that are arranged in a way defined by $\G$,
and the arithmetic group $\mathfrak G$ that
acts on $\sV_h$. Let $R_\G\to\sV_\G\times\sV_\G$
denote the groupoid induced from the action
of $\mathfrak G$
on $\sV_h$ and the closed embedding $\sV_\G\inj\sV_h$.
Via the Torelli theorem and the surjectivity
of the period map for K3 surfaces,
the points of the stack $\sV_\G/R_\G$ correspond to
$h$-tuples $(Y_1,...,Y_h)$ of Jacobian elliptic K3 surfaces,
one surface for each vertex of $\G$,
where adjacent surfaces (that is, surfaces
that correspond to adjacent vertices of $\Gamma$)
have isomorphic $\tD_4$-fibres.

Let $\sJE_\G$ denote the reduced closed substack of 
$\sJE^{sm,n.m.}_h$ whose geometric points are
configurations $x=\left( V_i\vert i\in\G\right)$.
So each $V_i$ is K3 and is smooth outside the double locus
$D_i=V_i\cap(\cap _{j\ne i}V_j)$. Each component of $D_i$ is
a fibre of type $\bD_4$. So $\sJE_\G$ equals
$(per^+)^{-1}(\sV_\G/R_\G)$.
Let $\sJE_{1^h}$ denote the substack 
defined by the condition that each $V_i$ is a special Kummer surface. 
So $\sJE_{1^h}$ equals $(per^+)^{-1}(\sW_{1^h}/R_{1^h})$,
where $R_{1^h}$ is the groupoid over $\sW_{1^h}$
induced from the the groupoid $R_\G$
(or, equivalently, from the action
of $\mathfrak G$ on $\sV_h$).
So there is a diagram with $2$-Cartesian squares
$$\xymatrix{
  {\sJE_{1^h}}\ar@{^{(}->}[d]\ar[r]^-{per^+_{1^h}}&{\sW_{1^h}/R_{1^h}}\ar@{^{(}->}[d]\\
  {\sJE_\G}\ar@{^{(}->}[d]\ar[r]^{per^+_\G}&{\sV_\G/R_\G}\ar@{^{(}->}[d]\\
  {\sJE^{sm,n.m.}_h}\ar[r]^-{per^+}&{\sV_h/\mathfrak G}
}$$
whose vertical arrows are closed embeddings.

\begin{remark} It is clear that the groupoid $R_{1^h}\to\sW_{1^h}\times\sW_{1^h}$
  is isomorphic to an action of
  the group $\mathfrak G_{1^h}:=(SL_2(\Z)\wr\Symm_g)\times SL_2(\Z)$.
\end{remark}

\begin{lemma} $\sJE_\G$ is smooth
along $\sJE_{1^h}$. 
\begin{proof} Recall first the easy fact that $\dim\sJE_\G=9h+9$.

Let $x=\left( V_i\vert i\in\G\right)$ be a point in 
the subvariety $\sJE_{1^h}$ of $\sJE_\G$.
Let $\tV_i\to V_i$ be the minimal resolution and $\tD_i\subset\tV_i$
denote the total transform of $D_i$. 
So $\tD_i$ is a sum of $\tD_4$ fibres. Let $\ts_i\subset\tV_i$
be the given section. Then the Zariski tangent space $T_x\sJE_\G$
is given by 
$$T_x\sJE_\G=\left\{(\xi_i)\in H^1(\tV_i,T_{\tV_i}(-\log(\tD_i+\ts_i))):
\xi_i\vert_{V_i\cap V_j}=\xi_j\vert_{V_i\cap V_j}\right\};$$
this is because each $H^1$ classifies first order deformations of $\tV_i$
that preserve the combinatorial structure of the configuration
$\tD_i+\ts_i$, and then we must impose the condition that
when the cross-ratio of the four marked points on each $\bD_4$
fibre varies, it does so in a way that is compatible with the fact that
it lies on $V_i$ and $V_j$.

Then, if $V_i$ is a vertex of $\G$ whose valency is $r$,
$$\dim H^1(\tV_i,T_{\tV_i}(-\log(\tD_i+\ts_i)))=20-(1+5r)+r-1=18-4r;$$
this is because the first Chern classes of the
curves in the configuration $\tD_i+\ts_i$ on $\tV_i$
are not linearly independent in $H^1(\tV_i,\Omega^1_{\tV_i})$,
but rather satisfy $r-1$ linear conditions.
So, if $\G$ has $\g_j$ vertices of valency $j$, then
$$\dim T_x\sJE_\G=\sum_1^4\g_j(18-4j)-(h-1)=9h+9,$$
since $\sum j\g_j=2e$, where $e=h-1$ is the number of edges in $\G$.
Therefore $\dim T_x\sJE_\G=\dim\sJE_\G$
and the smoothness is established.
\end{proof}
\end{lemma}

\begin{lemma} 
\part[i] 
$per^+_\G:\sJE_\G\to\sV_\G/R_\G$
is an isomorphism over a neighbourhood of
the generic point of $\sW_{1^h}/\mathfrak G_{1^h}$.
\part[ii] $\sV_\G$ is smooth along $\sW_{1^h}$.
\begin{proof} 
We show first
that the derivative $per^+_*$ of $per^+$ is injective
at all points $x$ of $\sJE_\G$.

As above, $x=\left( V_i\vert i\in\G\right)$ and 
$$T_x\sJE_\G=\left\{(\xi_i)\in H^1(\tV_i,T_{\tV_i}(-\log(\tD_i+\ts_i))):
\xi_i\vert_{V_i\cap V_j}=\xi_j\vert_{V_i\cap V_j}\right\}.$$
Choose a generator 
$\omega_i$ of $H^0(\tV_i,\Omega^2_{\tV_i})=H^0(V_i,\omega_{V_i})$;
then $per^+_*$ is the linear map defined by contraction
of $(\xi_i)$ against the various vectors $(0,...,0,\omega_i,0,...,0)$,
so is clearly injective.

It follows that $per^+_\G$ is {\'e}tale over a \nbd of
$\sW_{1^h}/\mathfrak G_{1^h}$
and then that $\sV_\G$ is smooth along $\sW_{1^h}$.
Finally, the surjectivity of the period map for K3 surfaces
completes the proof of the lemma.
\end{proof}
\end{lemma}

\begin{corollary}\label{Chak} (Chakiris \cite{C1}, \cite{C2})
    The generic Torelli theorem holds
    for simply connected
    Jacobian elliptic surfaces.
    \begin{proof} This is an immediate consequence
      of the fact that $per^+$
  is an isomorphism over a neighbourhood
  of $\sW_{1^h}$.
    \end{proof}
    \end{corollary}
  This proof of generic Torelli does not rely on
  Chakiris' version of a stable reduction theorem \cite{C1}
  but rather on Theorem \ref{10.x}, which derives from the
  MMP. However,
  beyond that, the proof relies on his ideas.

Define the vector bundle $E_\Gamma\to\sV_\Gamma$
by the property that its fibre over the point
$(Y_1,...,Y_h)$ of $\sV_\Gamma$
is the vector space
spanned by the $h-1$ matrices $\Pi_e$, each of rank one,
where $e=(i,j)$ runs over the edges of $\Gamma$ and,
in the notation of Proposition \ref{12.9}, 
$$\Pi_e=[\omega_{Y_i}(P_{ij}),-\omega_{Y_{j}}(P_{ji})]\otimes
[\underline{I}_i,\underline{I}_{j}].$$
This is a vector bundle of rank $h-1$.

\begin{theorem}\label{8.y}
\part[i] There is a branch $B_\Gamma$ of
$PL_h$ that contains $\sV_\Gamma$.

\part[ii]\label{defn of EGamma}
To first order $B_\Gamma$
is, in a neighbourhood of $\sW_{1^h}$, 
the vector bundle $E_\Gamma$.
\begin{proof}
Choose a smooth \nbd $\sV_\G^0$ of $\sW_{1^h}$ in $\sV_\G$
and let $E_\G^0$ denote the restriction of $E_\G$ to $\sV_\G^0$.
Then 
Proposition \ref{12.9} gives a family of surfaces of genus $h$
parametrized by $\sV_\Gamma^0\times S_{h-1}$,
where $S_{h-1}$ is an $h-1$-dimensional polydisc,
and the image of $\sV_\Gamma^0\times S_{h-1}$
under the period map equals, to first order,
precisely $E_\G^0$.
Since $\dim E_\Gamma=10h+8$, which is the number of moduli,
and, by Corollary \ref{Chak},
the period map is generically injective,
$E_\Gamma^0$ is, to first order, the image of some open piece of the moduli space.
\end{proof}
\end{theorem}
\begin{theorem}\label{main}
\part[i] The branch $B_\G$ of $PL_h$ 
is the unique branch of $PL_h$ that contains $\sV_\G$.

\part[ii]
To first order, the period locus $PL_h$ equals the union
$\cup_\G E_\G$ of the vector bundles $E_\G\to\sV_\G$
in a \nbd of $\sW_{1^h}$.
\begin{proof} Proposition \ref{11.a} shows that,
in particular, $per^+$ is proper over some \nbd $\sN_h$ of $\sW_{1^h}$
in $PL_h$. Set $\sJE^{0}_h=(per^+)^{-1}(\sN_h)$.
Note that $\sV_\G\subset PL_h$ and 
put $\sV_\Gamma^0=\sN_h\cap\sV_\G$.
(We could have used this choice of $\sV_\G^0$ in the proof of Theorem \ref{8.y}.)

Let $E_\G^0\to\sV_\G^0$ be the restriction of $E_\G$ to $\sV_\G^0$.
The plumbing construction of Section \ref{11} and the 
formula of Proposition \ref{11.8} for the derivative of the period matrix
show that $E_\G^0$ is, to first order, a closed subvariety of $\sN_h$.
That is,
there is, for each $\G$, a closed substack 
$\sF_\G$ of $\sJE^{0}_h$
such that $per^+$ induces an isomorphism $\sF_\G\to E_\G^0$
to first order.

Now on one hand
$\sJE^{sm,n.m.}_h$, and so $\sJE^{0}_h$, is smooth and,
on the other hand,
for each $x=(V_1,...,V_h)\in\sW_{1^h}$,
Theorem \ref{10.x} shows that
$(per^+)^{-1}(x)$ is a finite set which consists of 
exactly one point for each alkane $\G$ of genus $h$.
Therefore $(per^+)^{-1}(\sW_{1^h})\cong\coprod_\G\sW_{1^h}$
and, in a \nbd of $\sW_{1^h}$, 
$\sJE_h^{0}=\coprod_\G\sF_\G$.

So $PL_h=\cup_\G E_\G^0$ in a \nbd of $\sW_{1^h}$.
\end{proof}
\end{theorem}
\end{section}
\begin{section}{Fay's formulae for homologically trivial plumbings 
of curves}\label{Fay's formulae}\label{curves1}
For a homologically trivial Fay plumbing of curves, the arguments of
Section \ref{11} go through
to recover Fay's Corollary 3.2, as follows.

Suppose that 
$\sC\to\Delta_t$
is a homologically trivial 
Fay plumbing of $C_a$ to $C_b$ that identifies $a$ with $b$
and $(\omega^{(1)}(t),...,\omega^{(g)}(t))$
is a normalized basis of $H^0(\sC_t,\omega_{\sC_t})$. 

\begin{proposition}\label{9.1}
$$\omega^{(j)}(t) \equiv \omega^{(j)}_{C_i}+t\eta^{(j)}_{C_i}\pmod {t^2}$$
where 
\begin{enumerate}
\item $(\omega^{(j)}_{C_a})_{1\le j\le g_a}$ is
a normalized basis of $H^0(C_a,\omega_{C_a})$,
$(\omega^{(j)}_{C_b})_{g_a+1\le j\le g_a+g_b}$ is
a normalized basis of $H^0(C_b,\omega_{C_b})$
and $\omega^{(j)}_{C_i}=0$ otherwise,
and 
\item there is a unique element ${\tilde{\eta}}_i\in H^0(C_i,\omega_{C_i}(2i))$,
normalized by the requirements that, firstly,
$\int_{A_k}{\tilde{\eta}}_i=0$ for every $A$-cycle
$A_k$ on $C_i$ and that, secondly,
${\tilde{\eta}}_a={\quarter}(z_a^{-2}+\ \textrm{h.o.t.})dz_a$
while
${\tilde{\eta}}_b=-{\quarter}(z_b^{-2}+\ \textrm{h.o.t.})dz_b$,
such that there is an equality
$$[\eta_{C_i}]=
{\tilde{\eta}}_i
\left[
{\underline\omega}_{C_a}(a), 
-{\underline\omega}_{C_b}(b)
\right]$$
of row vectors of length $g$. The definition of these row vectors
is analogous to that of the vectors which are
defined immediately after Lemma \ref{up to}.
\end{enumerate}
\noproof
\end{proposition}
It follows that
$$\tau(\sC_t)\equiv
\left[
\begin{array}{cc}
{\tau(C_a)} & {0}\\
{0} & {\tau(C_b)}
\end{array}
\right]
+t\left[
{\underline\omega}_{C_a}(a), 
-{\underline\omega}_{C_b}(b)
\right]
\otimes{\underline{v}}
\pmod {t^2}$$
where ${\underline{v}}=[{\underline{v}}_a,{\underline{v}}_b]$ 
and ${\underline{v}}_i$
is the vector of integrals of the form
${\tilde{\eta}}_i$ around the $B$-cycles
on $C_i$.
Since, by the bilinear relations for integrals of the first kind, 
$\tau(\sC_t)$ is symmetric, it follows that 
$v=\lambda
\left[
{\underline\omega}_{C_a}(a), 
-{\underline\omega}_{C_b}(b)
\right]$ for some scalar $\lambda$ and
$$\tau(\sC_t)\equiv
\left[
\begin{array}{cc}
{\tau(C_a)} & {0}\\
{0} & {\tau(C_b)}
\end{array}
\right]
+\lambda t\left[
{\underline\omega}_{C_a}(a), 
-{\underline\omega}_{C_b}(b)
\right]^{\otimes 2}
\pmod{t^2}.$$
We can calculate $\lambda$ from the bilinear
relations for integrals of the second kind on $C_a$ with only poles at $a$: 
$$\sum_{l=1}^{g_a}\left(\int_{A_l}\phi\int_{B_l}\psi-\int_{A_l}\psi\int_{B_l}\phi\right)
=\int_{\gamma}f\psi$$ 
where $\g$ is a loop in $C_a$ around $a$
and $\phi=df$.
Taking $\psi=\omega^{(j)}_{C_a}$ and $\phi={\tilde\eta}_{a}$
gives, via the expansions above in terms of power series of $\omega_{C_a}^j$ 
and $\eta_{C_a}^j$,
$$\int_{B_j}{\tilde\eta}_{a}=\frac{\pi {\sqrt{-1}}}{2}\omega^{(j)}_{C_a}(a),$$
so that $\lambda=\frac{2\pi {\sqrt{-1}}}{4}.$
This differs from the formula given in Fay's Corollary 3.2, 
in which $\lambda={\quarter}$, 
because our $1$-forms are normalized
by the requirement that $\int_{A_l}\omega^{(j)}=\delta_{jl}$ and not
$2\pi {\sqrt{-1}}\delta_{jl}$.
\begin{remark}
As already mentioned, on p. 41 of \cite{F1} the minus sign 
in front of ${\underline\omega}_{C_b}(b)$ is missing.
The reason is a different choice of normalization
in the plumbing construction,
which replaces $z_b$ by $-z_b$.
Logically, however, there is no difference.
\end{remark}
\end{section}
\begin{section}
{Poincar{\'e}'s asymptotic period relations}\label{curves2}
Suppose that $E_1,...,E_g$ are disjoint
curves of genus $1$
and that $D$
is a copy of $\P^1$. Fix points $a_i\in E_i$ and a local co-ordinate $z_i$
on $E_i$ at $a_i$. On $D$ fix a point $\infty$ and a global
co-ordinate $u$ on $D-\{\infty\}$. Fix
distinct points $b_1,...,b_g\in D-\{\infty\}$ given
by $u-b_j=0$.

Then successively making Fay plumbings of the $E_i$ to $D$
using these data
in a way that identifies $a_j\in E_j$ to $b_j\in D-\{\infty\}$
leads to a family $\sC\to S$ of genus $g$ curves
over a $g$-dimensional polydisc $S=S_g=\Delta_{t_1,...,t_g}$
with co-ordinates $t_1,...,t_g$.
We fix a symplectic basis $(A_j,B_j)$ of each $H_1(E_j,\Z)$
and let $v_j$ be the corresponding normalized $1$-form on $E_j$.
We write $v_i(a_i)=\frac{v_i}{dz_i}(a_i)$;
this notation will be used many times.

The next result is due to Poincar{\'e} \cite{P} when $g=4$ and
to Fay \cite{F1} in general. According to Igusa (\cite{I}, p. 167)
Poincar{\'e}'s proof exploited the fact that the theta
divisor on a Jacobian is of translation type,
and so can not in any obvious way be extended to
an analysis of period matrices of varieties of higher dimensions.
However, Fay's approach, of which
we include some details
that he omitted, goes via his plumbing construction
and so can be extended. It turns out that hyperelliptic
curves are particularly interesting from this point of view.

As already mentioned, this result is derived in \cite{FGSM}
from their global results.

\begin{theorem}\label{asymptotic}
(``Poincar{\'e}'s asymptotic period relations'', \cite{F1}, p.45)
For $i\ne j$ the entries $\tau_{ij}$
of the period matrix $\tau$ of $\sC_t$ can be written as
$\tau_{ij}={\btau}_{ij}u_{ij}$
where $u_{ij}\equiv 1\pmod {(t)}$
and
$${\btau}_{ij}=\frac{2\pi {\sqrt{-1}}}{16}t_it_jv_i(a_i)v_j(a_j)/(b_i-b_j)^2.$$
\begin{proof} Induction on $g$.

Construct a family $\sC'\to S_1$ of curves of genus $1$ by plumbing $a_1\in E_1$
to $b_1\in D$. This has a normalized $1$-form
$\omega_1(\sC'_{t_1})$.
Then near a point of $D-\{b_1\}$ we have
$$\omega_1(\sC'_{t_1})={\quarter}t_1v_1(a_1)\frac{du}{(u-b_1)^2}+O(t_1^2).$$
Now construct a genus $2$ family $\sC\to\Delta_{t_1,t_2}$ by plumbing
$\sC'$ to $E_2$ in a way that identifies $b_2\in D$ with $a_2\in E_2$.
There are normalized $1$-forms $\Omega_j(\sC_{t_1,t_2})$ on $\sC_{t_1,t_2}$,
where $j=1,2$. Then, near a point in $E_2$, we have
$$\Omega_1(\sC_{t_1,t_2})= t_2\omega_1(\sC')(b_2)\eta_2+O(t_2^2),$$
where $\eta_2$ is a meromorphic $1$-form on $E_2$ with a double pole at $a_2$
and $\eta_2={\quarter}v_2(a_2)(z_2^{-2}+\ \textrm{h.o.t.})dz_2$.
Moreover, $\int_{A_2}\eta_2=0$.
So
\begin{eqnarray*}
\tau_{12}=\int_{B_2}\Omega_1(\sC)&\equiv&
t_2{\quarter}t_1v_1(a_1)\frac{1}{(b_2-b_1)^2}\int_{B_2}\eta_2\pmod{(t_1^2,t_2^2)}\\
&\equiv&\frac{2\pi{\sqrt{-1}}}{16}t_1t_2v_1(a_1)v_2(a_2)/(b_2-b_1)^2.
\end{eqnarray*}
Now assume that $g\ge 3$ and that the result is true for plumbings
of genus $\le g-1$. Write
$$\frac{2\pi {\sqrt{-1}}}{16}v_i(a_i)v_j(a_j)/(b_i-b_j)^2=c_{ij}.$$
Note that $c_{ij}\in\C^*$ since $v_i(a_i)\ne 0$.

Suppose that $k\in[3,g]$ and that $j\in[1,2]$. Then,
by induction,
$$\tau_{12}=c_{12}t_1t_2u_{12k}+t_kX_{12k}$$
where $u_{12k}\equiv 1\pmod{(t)}$ and $X_{12k}$
is some function.

Moreover, $\tau_{12}\equiv 0\pmod{t_j}$
since setting $t_j=0$ gives a plumbing where the curve $\sC_t$ remains singular
and one of its
irreducible components is the non-varying
curve $E_j$ of genus $1$.

Therefore $X_{12k}$ is divisible by $t_j$ and then we can write
$$\tau_{12}=t_1t_2(c_{12}u_{12k}+t_kY_{12k}).$$
Since $c_{12}\ne 0$ the induction is complete.
\end{proof}
\end{theorem}

It follows that
the off-diagonal quantities $y_{ij}={\btau}_{ij}^{-1/2}$ satisfy the
${g\choose 4}$ Pl{\"u}cker relations
$$y_{ij}y_{kl}-y_{ik}y_{jl}+y_{il}y_{jk}=0.$$
Clearing denominators gives an
octic polynomial $f_{ijkl}(\tau)$
for each Pl{\"u}cker relation, defined by
\begin{eqnarray*}
f_{ijkl}(\tau)&=&
2\tau_{ij}\tau_{kl}\tau_{il}\tau_{jk}\tau_{ik}\tau_{jl}
(\tau_{ik}\tau_{jl}+\tau_{il}\tau_{jk}+\tau_{ij}\tau_{kl})\\
&-&(\tau_{ij}^2\tau_{il}^2\tau_{jk}^2\tau_{jl}^2+
\tau_{ik}^2\tau_{il}^2\tau_{jk}^2\tau_{jl}^2+
\tau_{ij}^2\tau_{ik}^2\tau_{jl}^2\tau_{kl}^2),
\end{eqnarray*}
such that $f_{ijkl}(\btau)=0$.
Denote by $T$ the ideal generated by $\{\tau_{ij}\vert i\ne j\}$
and denote by $\sJ_g^{plumb}$ the locus in $\mathfrak H_g$ that consists
of the period matrices of the curves $\sC_t$ that arise
from the plumbing construction.
Then $f_{ijkl}(\tau)\equiv 0\pmod {T^9}$
on $\sJ_g^{plumb}$.  We say that ``the equations
$f_{ijkl}(\tau)=0$ 
are asymptotic period relations''.
That is, the functions $f_{ijkl}(\tau)$ reduce, modulo $T^9$,
to forms that vanish on the tangent
cone of $\sJ_g^{plumb}$ along $\sdiag_g$ inside $\mathfrak H_g$
Therefore there are holomorphic
functions $F_{ijkl}$ on $\mathfrak H_g$ such that each $F_{ijkl}$ 
vanishes on $\sJ_g^{plumb}$
and $F_{ijkl}$ is congruent to $f_{ijkl}$
modulo $T^9$.

We shall see later (Corollary \ref{plumb})
that, in a neighbourhood of the diagonal locus,
$\sJ_g^{plumb}$ coincides with the Jacobian locus.

Let $\sJ_g^c$
denote the closure
of the Jacobian locus $\sJ_g$ inside the stack $\sA_g$
along the closed substack $\sdiag_g$ of $\sA_g$
that parametrizes products of elliptic curves.
Note that $\sdiag_g$ is isomorphic to the quotient $(\sA_1)^g/\Symm_g$
of $\sA_1^g$ by the symmetric group. In transcendental terms,
$$\sA_g=\mathfrak H_g/Sp_{2g}(\Z)\ \textrm{and}\  
\sdiag_g=\mathfrak{Diag}_g/(SL_2(\Z)\wr\Symm_g).$$
Let $\bsM_g^c$ denote the open substack
of the stack $\bsM_g$ of stable curves of genus $g$
that parametrizes curves of compact type.
Then the Jacobian morphism $\sM_g\to\sA_g$ extends
to a proper morphism $\bsM_g^c\to\sA_g$ whose image
is $\sJ_g^c$.
We let $\mathfrak J_g$ and $\mathfrak J_g^c$ denote the
 inverse images in $\mathfrak H_g$
of $\sJ_g$ and $\sJ_g^c$.

We show next, in Proposition \ref{below}, that these ``asymptotic relations''
are exactly the defining equations, in terms of the entries
$\tau_{ij}$ of the period matrix, of the associated graded ring
belonging to the closed subvariety $\mathfrak{Diag}_g$ of $\mathfrak J_g^c$.

We shall use the term \emph{multi-elliptic}
to refer to a stable curve of genus $g$ that contains $g$
elliptic components (and maybe some smooth rational components).
Such a curve is necessarily of compact type.
Then $\sdiag_g$ is the image of the stack of multi-elliptic stable curves.

Let $\sG r_{\sJ_g^c}$ denote the sheaf of graded
$\sO_{\sdiag_g}$-algebras that is associated
to the closed embedding $\sdiag_g\inj\sJ_g^c$
and $\sG r_{\mathfrak J_g^c}$ its pullback to
$\sO_{\mathfrak{Diag}_g}$.

The Grassmannian $Grass(2,g)$ is embedded in 
$\Proj\C[\{y_{ij}\}]=\P^{{g\choose 2}-1}_y$
via the Pl{\"u}cker co-ordinates $y_{ij}$. Let $X$ denote the closure
of the image of $Grass(2,g)$ under the generically finite
rational map 
$$\P^{{g\choose 2}-1}_y-\to\P^{{g\choose 2}-1}_{\btau}$$
given by $y_{ij}\mapsto y_{ij}^{-2}=\btau_{ij}$ and let 
$\hX\subseteq\A^{g\choose 2}$ be the affine cone
over $X$ with affine co-ordinate ring $\C[\hX]$.
The quadratic Pl{\"u}cker identities relating
the $y_{ij}$ give rise to octic relations between the $\btau_{ij}$
that are the defining equations of $\hX$ (or of $X$).

\begin{proposition}\label{below}\label{12.2} $\sG r_{\mathfrak J_g^c}$ is isomorphic,
as a sheaf of graded $\sO_{\mathfrak{Diag}_g}$-algebras, to the
constant sheaf $\sO_{\mathfrak{Diag}_g}\otimes_{\C}\C[\hX]$.
\begin{proof}
Let $\sE\to\mathfrak{Diag}_g$ be the $g$-fold universal elliptic curve
over $\mathfrak{Diag}_g$,
so that $\dim\sE=2g$. Set $U=((\P^1)^g-\Delta)/PGL_2$
where $\Delta$ is the union of all the diagonals.
Let $S$ be a $g$-dimensional polydisc with co-ordinates
$t_1,...,t_g$ and put $L=\sE\times U\times S$. Then a suitable
Fay plumbing is a family of curves over $L$ which then gives a period map
$$h:L\to \mathfrak H_g$$
such that $\tau_{ij}={\btau}_{ij}u_{ij}$
as above.
Moreover, the sublocus of $L$ defined by $t_1=\cdots=t_g=0$
is everywhere locally the base of a miniversal
deformation of a singular stable curve $C_0$
which is of the form $C_0=\P^1+\sum_1^gE_i'$
for varying
curves $E_i'$ of genus $1$.
Therefore, by the local Torelli theorem for smooth curves,
there is a dense open subspace of $L$ along
which the derivative of $h$ is injective.
So the image of $h$ is open inside some branch of $\mathfrak J_g^c$
along $\mathfrak{Diag}_g$.
In particular, the dimension of $h(L)$ is $3g-3$.

\begin{lemma}\label{unibranch} $\sJ_g^c$ is unibranched
along $\sdiag_g$.
\begin{proof} Suppose that $E_1,...,E_g$ are curves of genus $1$
and that $\bsM_{\sum E_i}\subset\bsM_g$ is the locus of
multi-elliptic curves $C$ that contain each $E_i$.
That is, $C$ is a tree whose components are the curves $E_i$
together with some copies of $\P^1$.

To prove the lemma it is enough to show that $\bsM_{\sum E_i}$
is connected. We do this by induction on $g$.

This is clearly true when $g=1$. Say $\bsM'=\bsM_{E_g}$
and $\bsM''=\bsM_{\sum_1^{g-1}E_i}$. These are connected, by
the induction hypothesis.
Let $\bsM'_1$ and $\bsM''_1$ be their inverse images
in the stack $\bsM_{g,1}$ of $1$-pointed stable curves of genus $g$.
There are forgetful morphisms
$\bsM'_1\to\bsM'$ and $\bsM_1''\to\bsM''$.
Since their fibres are connected, both $\bsM'_1$ and $\bsM''_1$
are connected.
There is a clutching morphism
$\bsM'_1\times_{\C}\bsM''_1\to \bsM_{\sum_1^gE_i};$
since this is surjective the lemma is proved.
\end{proof}
\end{lemma}

\begin{corollary}\label{plumb} $\sJ_g^{plumb}$ coincides with $\sJ_g^c$ in
a neighbourhood of the diagonal locus.
\begin{proof} This follows from Lemma \ref{unibranch}
and the fact that $\sJ_g^{plumb}$ is exactly $h(L)$.
\end{proof}
\end{corollary}

Set $R=\sO_{\mathfrak H_g}$ and $S=\sG r_TR$, the sheaf
of graded $\sO_{\mathfrak{Diag}_g}$-algebras that
is associated to the ideal $T$.
Since the spaces concerned are Stein
we shall not emphasize the distinction between a ring and a sheaf of rings.
Given a graded $\sO_{\mathfrak{Diag}_g}$-algebra such as $S$,
$\Proj S$ will denote the projectivization relative to $\mathfrak{Diag}_g$
in the analytic category.

Let $I$ be the ideal of $R$ that defines $\mathfrak J_g^c$,
$J$ the ideal of $R$ generated by
the functions $F_{ijkl}$
and $K$ the ideal of $R$ that defines $(\mathfrak{Diag}_g\times \hX)\cap\mathfrak H_g$.
Note that $K$ is defined by the functions $f_{ijkl}$.
Let $\bI$, $\bJ$ $\bK$ be the associated graded ideals of $S$.

Note that $\bK$ is the defining ideal of $\mathfrak{Diag}_g\times X$
inside $\Proj S$ and so is prime.

By Corollary \ref{plumb}, $J\subset I$, so that $\bJ\subset\bI$.
Modulo $T^9$, the functions $F_{ijkl}$ and $f_{ijkl}$ are equal,
and so $\bJ=\bK$.
Therefore $\bK\subseteq\bI$.

Suppose that $\bK\ne\bI$.
Now $\Proj (S/\bK)=\mathfrak{Diag}_g\times X$,
which is reduced and irreducible of dimension $3g-4$.
So $\dim\Proj(S/\bI)\le 3g-5$.
On the other hand, $\Proj(S/\bI)$
is the exceptional divisor in the blow-up
$\Bl_{\mathfrak{Diag}_g}\mathfrak J_g^c$ and so has dimension
$3g-4$. Therefore $\bK=\bI$
and the proposition follows.
\end{proof}
\end{proposition}
\begin{corollary}
Up to graded equivalence
the singularity of $\sJ_g^c$
at a given point $x$ of $\sdiag_g$ is isomorphic to the cone $\hX$
and is independent of $x$.
\noproof
\end{corollary}
When $g=4$ (the case considered in detail by Poincar{\'e})
$\hX$ is then an octic hypersurface in $\A^6$ so that $\sJ_4^c$
does not have rational singularities along $\sdiag_4$.
\end{section}
\begin{section}{Hyperelliptic curves and alkanes}\label{curves3}
  We continue with the notation of the previous section,
  \emph{except} that ``elliptic curve'' will mean
  ``curve $E$ of genus $1$ provided with an involution $\iota$
  with $4$ fixed points''. By abuse of notation, we 
  denote $\iota$ by $[-1]_{E}$ and the fixed locus of $\iota$
  by $E[2]$.

Suppose instead that we choose
points $a_1$ on $E_1$, $b_{i-1}$ and  $a_{i}$ on $E_i$ for $i=2,...,g-1$
and $b_{g-1}$ on $E_g$
and a local co-ordinate
$z_x$ at each point $x$, and then plumb $E_i$ to $E_{i+1}$
by identifying $a_i$ to $b_i$
in a chain. This leads to a family $\sC\to S$
of genus $g$ curves where $S=\Delta_{t_1,...,t_{g-1}}$ is a
$(g-1)$-dimensional polydisc with co-ordinates
$t_1,...,t_{g-1}$.

Let $\tau_i$ be the period of $E_i$, defined, as before, after the choice of
a symplectic basis $(A_j,B_j)$ of $H_1(E_j,\Z)$ and normalized $1$-form
$\omega_j$ on $E_j$.
Construct a symplectic basis of $H_1(\sC_t,\Z)$
as a union of the symplectic bases $(A_j,B_j)$.

\begin{lemma} If each point $b_{i-1},a_i$ lies in $E_i[2]$
    and at each point $b_{i-1},a_i$
on $E_i$ a local co-ordinate is chosen that is anti-invariant
under the involution $[-1_{E_i}]$
then the family $\sC\to S$ is hyperelliptic.
\begin{proof} The constraint on the points
$a_1,...,a_{g-1}$ implies that the chain
$\sC_0=E_1\cup\cdots\cup E_g$ is a double cover
of a chain $\sB_0=\G_1\cup\cdots\cup\G_g$ of copies
of $\P^1$ and that the map $\sC_0\to\sB_0$
is not {\'e}tale at the nodes.
So, via Proposition \ref{plumbing involutions},
Fay plumbings can be constructed simultaneously
to give a double cover $\sC\to\sB$ over $S$ where $\sB\to S$
is a family of genus zero curves.
\end{proof}
\end{lemma}

We next give a description of
the closure $\sHyp_g^c$ of the hyperelliptic locus in $\bsM_g^c$
near $\sdiag_g$ in combinatorial terms.

\begin{lemma} Suppose that $C$ is multi-elliptic
  and that $E_1,...,E_g$ are its components of genus $1$.

\part[i] $C$ is hyperelliptic
if and only if
$C$ contains no rational curves,
each component $E_i$ of $C$ is elliptic
and contains at most $4$
nodes of $C$, and
each node of $C$ that lies on $E_i$
also lies in $E_i[2]$.

\part[ii] If $C$ is hyperelliptic then the hyperelliptic involution $\iota_C$
preserves each component $E_i$ of $C$ and
restricts to $[-1_{E_i}]$, and
the fixed point locus of $\iota_C$
consists of the nodes of $C$ and
each $E_i[2]$.
\begin{proof}
This is an elementary exercise.
\end{proof}
\end{lemma}
So we
can describe stable multi-elliptic hyperelliptic curves $C$
in terms of alkanes of the same genus, where the genus of an alkane
is the number of carbon atoms in it.
\begin{corollary}
Each fibre in $\sHyp_g^c$ over a point in $\sdiag_g$
is a finite set that corresponds naturally to the set of alkanes of genus $g$.
\begin{proof}
This is a translation of the preceding lemma.
The carbon atoms correspond to the elliptic curves $E_i$,
the bonds on the atom correspond to the points in $E_i[2]$
and the hydrogen atoms to the points in $E_i[2]$
that are not nodes.
That is, the the hydrogen atoms correspond to the
fixed points of the hyperelliptic involution
that are not nodes.
\end{proof}
\end{corollary}

\begin{corollary} The branches of $\sHyp^c_g$ through $\sdiag_g$
in $\sA_g$ correspond to the alkanes of genus $g$.
\noproof
\end{corollary}

We shall also use the term alkane to refer to a hyperelliptic
multi-elliptic stable curve.

Suppose that $C=\sum_1^gE_i$ is an alkane.
Suppose that $K$ is the set of edges,
so that $\{i,j\}\in K$ if and only if $i\ne j$ and $E_i$ meets $E_j$.
Choose a local co-ordinate
on each $E_i$ at each node of $C$
(not necessarily anti-invariant under the involution $[-1_{E_i}]$)
and let $\sC\to S_{g-1}=\Delta_{\{t_k\vert k\in K\}}$
be the result of plumbing together the curves
$E_i$ according to these data.

\begin{theorem} Modulo $(\{t_l\vert l\in K\})^2$,
the off-diagonal entry
$\tau(\sC_t)_k$ of the period matrix of $\sC_t$
is a multiple of $t_k$
if $k\in K$. The other off-diagonal entries vanish.
\begin{proof} We can suppose the curves $E_i$ ordered so that
$E_g$ meets only $E_{g-1}$
and then argue by induction on $g$.

Suppose that $a\in E_{g-1}$
is identified with $b\in E_g$.
Suppose that $\sC'\to S'=\Delta_{t_1,...,t_{g-2}}$
is the result of plumbing $E_1,...,E_{g-1}$
according to the data,
so that $\sC\to S$ is the plumbing of $\sC$ to $E_g\times S'$
that identifies $\{a\}\times S'$ with $\{b\}\times S'$.
Write $(t')=(t_1,...,t_{g-2})$ and $(t)=(t',t_{g-1})$ and suppose that
$(\Omega_1',...,\Omega_{g-1}')$ is a basis of
normalized $1$-forms on $\sC'$. By Fay's formula,
$\tau(\sC_t)$ is congruent to
$$\left[
\begin{array}{cc}
{\tau(\sC'_{t'})}&{0}\\
{0}&{\tau_{g}}
\end{array}
\right]
+{\frac{\pi {\sqrt{-1}}}{2}}t_{g-1}[\Omega_1'(\{a\}\times S'),...,
\Omega_{g-1}'(\{a\}\times S'),
-\omega_g(b)]^{\otimes 2}$$
modulo $(t)^2$.
The restriction $\Omega_i'\vert_{E_{g-1}}$
is identically zero for all $i\le g-2$,
and so $\Omega'_i(\{a\}\times S')\equiv 0\pmod{t'}$ for $i\le g-2$.
Moreover, according to the induction hypothesis,
$\tau(\sC'_{t'})$ has the required form,
and the result follows.
\end{proof}
\end{theorem}

Let $\mathfrak{Hyp}_g^c$ denote the closure of the
hyperelliptic locus in $\mathfrak H_g$.

\begin{theorem}\label{4.10} (Asymptotic period
relations for the hyperelliptic locus)
To first order in a neighbourhood of $\mathfrak{Diag}_g/\Symm_g$
each branch of $\mathfrak{Hyp}_g^c$ is described as a subvariety
of $\mathfrak H_g$
by the vanishing of the entries $\tau_k$ in the period matrix $\tau$
where $k$ runs over the set of pairs that are not edges of the corresponding alkane.

In particular,
the branch corresponding to the linear alkane
equals, to first order, the locus $\mathfrak T_g$ of tridiagonal matrices.
\begin{proof}
We prove this only for the linear alkane.
The proof in general is the same but with more
complicated notation.

Let $E_\tau=\C/(\Z\tau\oplus\Z)$ be the elliptic curve
corresponding to $\tau\in\mathfrak H_1$.
Pick a co-ordinate $z$ on $\C$ such that the involution
$[-1_{E_{\tau}}]$ acts via $z\mapsto -z$.
Let $\sF\to\mathfrak{Diag}_g$ be the family whose fibre over
$(\tau_1,...,\tau_g)$ is $E_{\tau_2}\times\cdots\times E_{\tau_{g-1}}$.
Then the Fay plumbing just described, using the
co-ordinates provided, gives a family
$\sC\to\sF\times S_{g-1}$ of stable curves which is minimally
versal at each point where $t_1=\cdots=t_{g-1}=0$
and so is minimally versal everywhere.

Fixing a non-zero $2$-torsion point on each of $E_2,...,E_{g-1}$
defines a section of $\sF\to\mathfrak{Diag}_g$.
Then the restriction
of $\sC$ to the corresponding subvariety $\mathfrak{Diag}_g\times S_{g-1}$
of $\sF\times S_{g-1}$
is then, over the complement of the discriminant,
a family of hyperelliptic curves.

Since $\dim(\mathfrak{Diag}_g\times S_{g-1})=2g-1$ this is
everywhere a minimally versal family of hyperelliptic curves.
Then, by the local Torelli theorem for hyperelliptic curves,
the image $T$ of $\mathfrak{Diag}_g\times S_{g-1}$ in $\mathfrak H_g$
is of dimension $2g-1$. To first order $T$ lies in $\mathfrak T_g$,
and the theorem is proved.
\end{proof}
\end{theorem}

\end{section}

\bibliography{alggeom,ekedahl}
\bibliographystyle{pretex}
\end{document}